\documentclass[12pt]{amsart}

\usepackage{amsfonts}
\usepackage{epsf}
\usepackage{epic}
\usepackage{eepic}
\usepackage{latexsym}

\def\lmn#1{\vadjust{\setbox1=\vtop{\hsize 12mm
\parindent=0pt\baselineskip=9pt
\rightskip=4mm plus 4mm#1}
\hbox{\kern-12mm\smash{\raise .5ex\box1}}}}

\newtheorem{theorem}{Theorem}[section]

\newtheorem{lemma}[theorem]{Lemma}
\newtheorem{proposition}[theorem]{Proposition}
\newtheorem{definition}[theorem]{Definition}
\newtheorem{corollary}[theorem]{Corollary}
\newtheorem{remark}[theorem]{Remark}
\newtheorem{remarks}[theorem]{Remarks}

\def\C{{\mathbb C}}
\def\L{{\mathcal{L}}}
\def\S{\Sigma}
\def\Z{{\mathbb Z}}
\def\a{{\alpha}}
\def\b{{\beta}}
\def\g{{\gamma}}

\def\E{{\mathcal{E}}}
\def\U{{\mathcal{U}}}
\def\G{{\mathcal{G}}}
\def\J2{J^{(2)}}

\def\cO{{\mathcal{O}}}

\def\bP{{\mathbb P}}

\def\ra{\rightarrow}
\def\P{{\mathcal P}}
\def\M{{\mathcal M}}
\def\F{{\mathcal F}}
\def\G{{\mathcal G}}

\def\pia{{\pi_{\a*}}}
\def\pib{{\pi_{\b*}}}
\def\pig{{\pi_{\g*}}}

\def\lag{{\langle}}
\def\rag{{\rangle}}

\newcommand{\tr}{\mathop{\fam0 Tr}\nolimits}
\newcommand{\ch}{\mathop{\fam0 Ch}\nolimits}
\newcommand{\td}{\mathop{\fam0 Td}\nolimits}
\newcommand{\Pic}{\mathop{\fam0 Pic}\nolimits}
\newcommand{\Gr}{\mathop{\fam0 Gr}\nolimits}
\newcommand{\Nm}{\mathop{\fam0 Nm}\nolimits}
\newcommand{\Arf}{\mathop{\fam0 Arf}\nolimits}

\def\noi{\noindent}
\def\v8{\vskip 8pt}

\def\build#1_#2^#3{\mathrel{\mathop{\kern 0pt#1}\limits_{#2}^{#3}}}
\def\mapright#1{\smash{\mathop{\longrightarrow}\limits^{#1}}}

\newcommand{\Div}{\mathop{\fam0 Div}\nolimits}
\newcommand{\im}{\mathop{\fam0 Im}\nolimits}

\newcommand{\Jt}{J^{(2)}}
\def\Jf{J^{(4)}}
\def\Jr{J^{(r)}}

\def\tS{\widetilde{\S}}

\def\ch{{\widetilde {Ch}}}
\def\pa{{p_\a}}
\def\pb{{p_\b}}
\def\pg{{p_\g}}

\begin{document}

\title[Involutions on Moduli Spaces]{Involutions on Moduli Spaces and 
Refinements of the Verlinde Formula}

\author{J{\o}rgen Ellegaard Andersen}
\address{Dept. of Math., Univ. of Aarhus, DK-8000 Aarhus C, Denmark}

\email{andersen@mi.aau.dk}

\author{Gregor Masbaum}
\address{Institut de Math\'ematiques de Jussieu, Equipe `Th\'eories 
G\'eom\'etriques', Case 7012, Universit\'e Paris VII, 75251 Paris Cedex 05, 
France}
\email{masbaum@math.jussieu.fr}
\date{November 2, 1998 (First Version October 22, 1997)}

\thanks{Research at MSRI is supported in part by NSF grant DMS-9022140.}
\thanks{J.E.A. is supported in part by NSF grant DMS-93-09653 and by the
Danish Research Council.}

\begin{abstract}  The moduli space $M$ of semi-stable rank $2$ bundles
  with
 trivial
  determinant over
  a complex curve $\S$ carries involutions naturally
  associated 
to $2$-torsion points on the Jacobian
  of the curve. For every lift of a $2$-torsion point to  a
  $4$-torsion point,
  we 
define a  lift  of the  involution to the
  determinant line bundle $\L$. We obtain an explicit
  presentation of the group generated by these lifts  in terms of
  the 
order $4$ Weil
  pairing. This is related to the triple intersections of  the components
  of the fixed point
  sets in $M$, which we  also determine completely using the  order
  $4$ Weil
  pairing. The lifted involutions  act on the spaces 
  of holomorphic sections of powers of $\L$, whose dimensions are
  given by the
  Verlinde formula. We compute the characters of
  these  vector spaces as representations of the group generated by
  our lifts, and we obtain  an explicit isomorphism 
(as group representations)  with
  the combinatorial-topological TQFT-vector spaces of \cite{BHMV}. As
  an application, we describe 
  a `brick decomposition', with explicit dimension formulas, of the
  Verlinde vector spaces.  We also obtain similar results in the
  twisted ({\em i.e.}, degree one) case. 
\end{abstract}

\maketitle

\tableofcontents

\section*{Introduction and Motivation.} The celebrated Verlinde
formula
 \cite{Ve} 
gives the
dimension of certain vector spaces of so-called `conformal blocks' appearing
  in conformal field theory.  In this paper, we will take the point of
  view of algebraic geometry and think of the conformal blocks
  as
  holomorphic sections of powers of the determinant line bundle over
  moduli spaces of  
semi-stable bundles  with fixed determinant over a complete
non-singular complex curve
$\S$. According to Atiyah
  \cite{At} and Witten
\cite{Wi}, the spaces of holomorphic sections should also fit into
 $2+1$-dimensional
`Topological Quantum Field Theories' (TQFT). This geometric construction has 
been studied quite a lot (see
{\em e.g.} \cite{ADW,MS,RSW,Hi,CLM,Th,Th2,Sz,BSz,BL,KNR,Fa}). 

In \cite{BHMV},  a combinatorial-topological construction of
TQFT-functors was described, based on a particularly simple
construction of 
the
Witten-Reshetikhin-Turaev $3$-manifold  invariant \cite{RT} in the
$SU(2)$-case through
the  Kauffman bracket. In that paper one also constructed certain  
involutions on the TQFT-vector spaces which were then used to 
decompose the 
vector spaces into direct summands. These involutions are associated
to
 simple closed
curves on the underlying smooth surface of $\S$, and they generate a
kind 
of Heisenberg group
presented in terms of the $mod $ $4$ intersection form.  These ideas
were 
 developed
further  in \cite{BM} to construct
spin-refined TQFT's. 

The starting point  for the present paper was the idea that the involutions of
\cite{BHMV} should correspond on the algebraic-geometric side to the 
 involutions on moduli space naturally associated to $2$-torsion
points on the Jacobian of $\S$. More precisely, the $2$-torsion
points  should correspond  to the $mod $ $2$  homology classes of the
simple
 closed
curves. Note that  the involutions on the spaces of holomorphic
sections require a  choice of lift to the determinant line bundle
$\L$. It is easy to see that these lifts generate a central 
extension of the group of $2$-torsion
points whose  alternating form  is given
by the order $2$ Weil pairing. It follows that this extension is
indeed 
abstractly isomorphic
to the one that appeared in
\cite{BHMV}. The ambiguity in the choice of lift is reflected on the
topological side by the choice of a simple closed curve within its $mod
$ $2$  homology class.

One of our motivations in this paper is to establish  a more precise  
correspondance
between the two theories. The key idea is to define, for every lift of
a $2$-torsion
point $\a$ to a  
$4$-torsion point $a$, an involution $\rho_a$ on $\L$ which covers the
action of $\a$  on the moduli space $M$. The sign of this lift
$\rho_a$ is fixed  by requiring it to act as the
identity over a certain component of the fixed point set of the involution
acting on $M$; this component is simply the one containing the class
of the semi-stable bundle $L_a\oplus L_a^{-1}$ (see section
\ref{prel} for more details). We will see in section
\ref{taugamma} that the action of this
lift 
  $\rho_a$ on holomorphic sections corresponds precisely to 
the involution $\tau_\g$ which is
associated in the \cite{BHMV}-theory to a simple closed curve $\g$
whose homology class is  Poincar\'e dual to  $a$. In this way, we
obtain  an
 explicit isomorphism 
(as group representations)  with
  the  TQFT-vector spaces of \cite{BHMV}. We believe this constitutes a
nice confirmation that there should be a natural  correspondence
between
 the two theories. 

 The main work in this paper is, however, algebraic-geometric in
  nature and consists of a detailed study of our lifts $\rho_a$ and
  their action on the Verlinde vector spaces. Our main results are as
  follows (see section \ref{smr} for more
  complete statements).   We show in theorem \ref{1.1} that our lifts 
satisfy $$\rho_a \,\rho_b\,=\,\lambda_4(a,b) \,\rho_{a+b}, $$ where 
$\lambda_4$ is the order $4$ Weil
  pairing (which is the algebraic-geometric analogue of the {\em mod}
  $4$ intersection form).  Along the way, we also determine in theorem
  \ref{evencasei}   the
  triple intersections of the components of the fixed point sets of
  the 
action on the various
  moduli spaces. The {\em r\^ole} played  by the order $4$ Weil
  pairing in this context doesn't seem to have been observed before.
We then compute in theorem \ref{1.2} the trace of the induced
  involutions  
$\rho_a^{\otimes
  k}$ on the spaces of holomorphic sections of the $k$-th
tensor power of $\L$. This determines the characters of these vector
  spaces as representations of the group generated by our lifts. We
  also 
obtain similar results in the
  twisted ({\em i.e.}, degree one) case, where the situation is somewhat
  simpler, as only the order $2$ Weil pairing is needed.

A corollary of our results
  is a `brick decomposition' of the
  spaces of holomorphic sections, the structure of which depends on
  the value of the level $k$ modulo $4$. This is, of course, analogous to 
the decomposition in
  \cite{BHMV}, but we will derive it, as well as explicit
  dimension formulas,  directly from the
  character of the representation.  At low levels, similar
  decompositions  have appeared  previously in the work of Beauville
  \cite{Be2} (see also Laszlo \cite{L}) for $k=2$,   and of  van
  Geemen 
and Previato \cite{vGP1,vGP2}, Oxbury and Pauly
\cite{OP}, and Pauly \cite{P} for $k=4$ (in our notation). Their
  approach is based
  on the 
relationship
  with abelian theta-functions and seems  quite  different
  from ours.  

This paper is organized as follows. After giving the basic
definitions in section \ref{prel}, we state our main results in
section \ref{smr}.  The relationship with the \cite{BHMV}-theory is
discussed in more detail in section \ref{taugamma}, and the  `brick 
decompositions' are
described in section
  \ref{3}.  The
  remainder of the paper is devoted to the proofs. (See remark
  \ref{logstruct} for the logical structure of the proofs.) The reader
  interested only in the algebraic geometry may
  skip section \ref{taugamma}, and no familiarity with 
\cite{BHMV} is necessary to understand the results and their proofs. 
{\small \v8\noi{\bf Acknowledgment.} 
We would like to thank Ch. Sorger for discussions on this project.}

\section{Basic definitions and notation.}\label{prel}

Let $\Sigma$ be a complete, non-singular curve over the complex
numbers 
of genus
$g\geq 2$. Let $M_d$ be the moduli space of
semi-stable bundles  of rank $2$ and degree $d\in \Z$ on $\Sigma$. 
There is a natural
algebraic action of the degree zero Picard group
$\Pic_0(\S)$ on $M_d$ gotten by tensoring. In this paper, we will use
the standard 
identification of  $\Pic_0(\S)$ with the Jacobian
$J(\S)$ and speak of an action of $J(\S)$ on $M_d$.

Fix a point
$p\in \Sigma$ and let $[p]$ be the associated line bundle. We
have the determinant morphism $\det : M_d \rightarrow \Pic_d(\S)$. We
put $M=\det^{-1}(\cO_\S)$ and $M'=\det^{-1}([p])$. 

As it is explained in \cite{DN} there is a
natural construction of a determinant line bundle $\L_d$ over
$M_d$. We will denote the restriction of
$\L_0$ to
$M$ by $\L$ and the restriction of $\L_1$ to
$M'$ by $\L'$. These two determinant bundles are generators
respectively of the Picard group of $M$ and of $M'$ (see
\cite{DN}).

 Let $k\geq 1$ be an integer called the {\em level}. Put 
$$Z_k(\Sigma) = H^0(M,\L^k)$$ 
$$Z'_k(\Sigma) = \left\{\begin{array}{ccl} 
H^0(M',{\L'}^{k/2}) & \ {\rm if} \ k\equiv 0 \ {\rm mod} \ 2\\
0& \ {\rm if} \ k\equiv 1 \ {\rm mod} \ 2\end{array}\right.$$

It is by now well-known that the dimensions of these vector spaces are
given by the
celebrated Verlinde formulas (to be recalled in section \ref{3}). 
Following   Thaddeus \cite{Th}, we will refer to $Z'_k(\S)$
as the {\em twisted case}, and to $Z_k(\S)$ as the {\em untwisted
  case.}

Let $\Jr$ be the subgroup of order $r$ points on $J(\S)=\Pic_0(\S)$. We denote by
 $L_\a$ (resp. $L_a$) the line bundle on
$\S$ corresponding to $\a\in
\Jt$ (resp. $a\in \Jf$). Note that the group  $\Jr$ is identified 
with $H^1(\Sigma;\mu_r)$, where $\mu_r\subset \C$ is the group of
$r$-th roots of unity.

The action of $\a\in  \Jt$ on $M_d$  is given by tensoring with
$L_\a$. Since $L_\a^{\otimes 2}\cong \cO_\S$, this preserves $M$ and
$M'$.  By abuse of notation, the automorphisms of $M$ and
$M'$ induced by $\a\in\Jt$ will again be denoted by $\a$.

\v8\noi{\bf The lifts $\rho_a$ and $\rho'_\a$.} By a lift of $\a$ to
$\L$  we mean an invertible bundle map from  $\L$ to itself covering
$\a$. For a lift to exist it suffices that $\a^*\L\cong \L$. This
is the case for every  $\a\in
\Jt$, since  the Picard
group of $M$ is  isomorphic to $\Z$, and $\L$ is ample, so that $\a$
must act trivially on the Picard group.  Therefore the
action of any
$\a\in
\Jt$ can be lifted to $\L$, and also to $\L'$, for the same
reason. Since 
the only algebraic 
functions on $M$ and $M'$ are the constant ones, it
follows that we can actually choose involutive lifts of each element of
$\Jt$. Any two involutive lifts of $\a$  agree up to sign.

To fix the signs of the lifts of $\a$ to $\L$ and $\L'$, we use the
fact that the sign can be read off in the fiber over a fixed point of
$\a$. Given $0\neq \a\in \Jt$, we 
use the notation $|X|_\alpha$ for the fixed point variety of the
automorphism induced by 
$\a$  on the various moduli spaces $X$. 

It is well-known that $|M'|_\a$ is isomorphic to the {\em Prym variety
  $P_\a$} associated to $\a$, and that $|M|_\a$ is isomorphic to the
  disjoint union of two
  copies of the {\em Kummer variety} $P_\a/\lag \pm 1\rag$. This result is
  due to Narasimhan and Ramanan \cite{NR}. In particular, $|M'|_\a$ is
connected and
non-empty, and  $|M|_\a$ has two
components. For any $a\in J^{(4)}$ such that $2a=\a$, we denote by
  $|M|_a^+$ 
the component of $|M|_\a$ containing
the S-equivalence class of the semistable bundle $L_a\oplus
L^{-1}_a$. (Note that this is indeed fixed under tensoring with
$L_\a$, since $L_a\otimes L_\a\cong L_a^{-1}$.) The other component
of $|M|_\alpha$ is denoted by $|M|_a^-$.

\begin{definition}\label{def11}
For $0\neq \a\in J^{(2)}$, we define $\rho'_\alpha$ to be the involutive lift  
to $\L'$ of $\alpha$ acting 
on $M'$  such that $\rho'_\alpha$ acts by {\em minus} the identity on
the restriction of $\L'$
to $|M'|_\alpha$.
\end{definition}

\begin{definition}\label{def12} For $a\in J^{(4)}$ such that $2a=\a\neq 0$, we define
$\rho_a$ to be the involutive lift  to $\L$ of $\alpha$ acting 
on $M$  such that  
$\rho_a$ acts by the identity on the restriction of $\L$
to the component $|M|^+_a$ of $|M|_\alpha$ specified by  $a$. 
\end{definition}
It will be convenient to extend this definition by letting
$\rho'_\alpha$ (resp. $\rho_a$) be the identity if $\a=0$
(resp. $2a=0$).

 \v8\noi{\bf Note.} According to Theorem F in \cite{DN}, we have that 
${\L'}^{-2} \cong K$, 
where $K$
is the canonical bundle of $M'$. The natural action of $\a$ on $K$
coincides with the one induced by both possible involutive lifts $\pm
\rho'_\alpha$.  This is because $\a$ acts as the identity on the fiber
of $K$ 
over the fixed point set $|M'|_\alpha$, since  $|M'|_\alpha$ has even
codimension (see section \ref{fpv}).  Similar comments apply in
the untwisted case. 
\v8
\noi{\bf The Weil pairing.}  (See {\em e.g.} \cite{Ho}.\footnote{We thank
  A. Beauville for pointing out this reference.}) Let $\M(\S)$ be the
field of meromorphic
functions on $\S$. The divisor of $f\in \M(\S)$ is denoted by
$(f)$. As already mentioned, we consider  $J^{(r)}$ to
be the  $r$-torsion points on the group
$\Pic_0(\S)$ which is naturally identified with
$\Div_0(\S)/\Div^{pr}(\S)$. (Here, $\Div_d(\S)$ is the group of degree
$d$ divisors, and
$\Div^{pr}(\S)$ are the principal divisors.) If $D=\sum n_j x_j$ is a 
divisor and $f\in
  \M(\S)$ a meromorphic function, we put  $f(D)=\prod
  f(x_j)^{n_j}$. The Weil pairing 
$$\lambda_r :
J^{(r)}\times J^{(r)} \ra \mu_r$$ is defined as follows. Given $a,b\in
J^{(r)}$, represent them  by divisors
$D_a,D_b$ with disjoint support, and pick $f,g\in \M(\S)$ such that
$(f)=r D_a$ and $(g)=r D_b$. Then
\begin{equation}
\lambda_r(a,b)=\frac{g(D_a)}{f(D_b)}.\nonumber 
\end{equation} The Weil pairing is antisymmetric and
non-degenerate. The fact that it takes values in $\mu_r$ follows from
Weil reciprocity (see \cite{GH}, p. 242).

\section{Statement of the main results.}\label{smr}

Let $\G(\J2,\L)$  be the group of automorphisms
of the  determinant line bundle $\L$  covering the action of $\Jt$ on
$M$. Since $\alpha^*(\L)\approx \L$ for every $\alpha\in
\Jt$, the group $\G(\J2,\L)$ is a central extension 
\begin{equation} \C^* \ \longrightarrow \  \G(\J2,\L)\
  \longrightarrow \ \Jt \label{tildE}. 
\end{equation} The same holds for the group $\G(\J2,\L')$ of automorphisms
of $\L'$ covering the action of $\Jt$ on
$M'$.  

\v8\noi{\bf Notation.} Let $\E\subset\G(\J2,\L)$ be  the subgroup generated by the
involutions $\rho_a$ ($a \in \Jf$), and let $\E'\subset\G(\J2,\L')$ be
the subgroup generated  by the involutions $\rho'_\alpha$ ($\a\in \Jt$). (See
definitions \ref{def12} and \ref{def11}.)
\v8

Our first main result gives a presentation of the groups $\E$ and
$\E'$, as follows.

\begin{theorem}\label{1.1}
The involutions $\rho_a$ and $\rho'_\alpha$ satisfy the following
relations: 
\begin{equation} \rho_a \,\rho_b\,=\, \lambda_4(a,b)\, \rho_{a+b}\label{rhoa}
\end{equation}
\begin{equation} \rho'_\alpha\,
\rho'_\beta\,=\, \lambda_2(\a,\b)\, \rho'_{\alpha+\beta}\label{rhoprima} 
\end{equation}
\end{theorem}

It follows that the group $\E$ is a central extension 
\begin{equation} \mu_4 \ \longrightarrow \  \E \ \longrightarrow \
  \Jt.\label{tildE2}
\end{equation}
This extension is non-trivial, since the associated alternating form  on $\Jt$ is the order $2$
Weil pairing $\lambda_2$. Indeed, this form  is
given by the commutator 
$$ c(\a,\b)=\rho_a \rho_b\rho_a^{-1}
\rho_b^{-1}=(\rho_a \rho_b)^2=\lambda_4(a,b)^2=\lambda_2(\a,\b).$$ The group $\E'$ is a
trivial extension of $J^{(2)}$, {\em i.e.}, it is isomorphic to $\mu_2
\times J^{(2)}$. 

\begin{remarks} {\em (i) The group $\G(\J2,\L)$ is known as the {\em
      Heisenberg group}. The fact that its alternating form is the order $2$
Weil pairing $\lambda_2$ is well-known. For example, it follows already from Beauville's isomorphism of $Z_1(\S)$ with the space of abelian theta-functions 
\cite{Be0}. This fact is also very easy to see from our point of view
(see remark \ref{coc}). Note that the alternating form determines the extension
(\ref{tildE}) (but not the extension (\ref{tildE2})) up to
isomorphism. 

(ii) If the alternating form is known, one knows {\em a priori} that
      our  lifts $\rho_a$ satisfy $\rho_a\rho_b=\pm
      \rho_{a+b}$ if $\lambda_2(\a,\b)=1$, and $\rho_a\rho_b=\pm i 
      \rho_{a+b}$ if
$\lambda_2(\a,\b)=-1$, where $2a=\a$ and $2b=\b$.  (This is because
      the lifts $\rho_a,\rho_b,\rho_{a+b}$ are involutions.) But the sign
      of the prefactors $\pm 1$ and $\pm i$ in these relations is, of
      course, not determined by the
      alternating form. The   contribution of theorem
\ref{1.1} is to show that with our `geometric' choice of lifts in terms of
      components of the fixed point sets, the
prefactors $\pm 1$ and $\pm i$ are given by the order $4$ Weil
pairing. Similar remarks apply in the twisted case.
}\end{remarks}

\noi{\bf Note.}  In the literature, the Heisenberg
group $\G(\J2,\L)$ is often described explicitly in terms of a `theta-structure' (see {\em
  e.g.} \cite{Be2}), that is, $\G(\J2,\L)$ is described as a
certain group structure on the set $\C^*\times (\Z/2)^g\times
(\Z/2)^g$. A theta-structure allows one to write  the extensions
(\ref{tildE}) and (\ref{tildE2}) as push-outs of an extension of $\Jt$ by
$\mu_2$. But this `reduction' to an extension by $\mu_2$ is not
canonical, as it depends on the choice of theta-structure (which comes down,
essentially, to the  choice  of a symplectic basis of $\Jt\approx
H^1(\S,\Z/2)$).\footnote{The {\em existence} of such a `reduction'
  follows already from the fact that the alternating
  form $c(\a,\b)=\lambda_2(\a,\b)$ takes values in $\mu_2$.  But there are many choices for
  this `reduction'.} {}From our point of view, such a choice is
neither necessary nor useful, as it  would break the symmetry of our
description of the group
$\E$, which is defined completely intrinsically in terms of the
involutions 
$\rho_a$. Therefore we  won't use theta-structures in this paper.

\v8

Theorem \ref{1.1} is related to the triple intersections  of the components of
the fixed point set on $M$. In fact, a key step in the proof is the
following result which does not seem to have
been observed before. 

\begin{theorem}\label{evencasei} Assume that $\a$ and $\b$ are
  distinct non-zero elements of $J^{(2)}$  such
that $\lambda_2(\a,\b)=1\in\mu_2$.  Let $a,b\in J^{(4)}$ such that $2a=\a$ and
$2b=\b$. Then 
$$ |M|^+_a \cap |M|^+_b \cap |M|^+_{a + b} \neq \emptyset \
\Leftrightarrow\  \lambda_4(a,b)=1\in\mu_4$$
$$ |M|^+_a \cap |M|^+_b \cap |M|^-_{a + b} \neq \emptyset \
\Leftrightarrow\  \lambda_4(a,b)=-1\in\mu_4$$
\end{theorem}

\noi{\bf Note.} Given theorem \ref{evencasei},  theorem \ref{1.1} in
the case $\lambda_4(a,b)=1$
follows immediately. Indeed, if the
triple intersection $|M|_a^+\cap|M|_b^+\cap |M|_{a+b}^+$ is non-empty,
it follows  from the definition of the lifts that $\rho_{a+b}=\rho_a
\rho_b$ (since $\rho_a \rho_b 
\rho_{a+b}$ must be a constant, and this constant can be computed in
the fiber of $\L$ over a triple intersection point.) The proof in
the general case, however, requires some further arguments. The most
interesting case is when $\lambda_4(a,b)=\pm i$. In this situation,
the fixed point sets $|M|_\a$, $|M|_\b$, $|M|_{\a+\b}$ don't
intersect, but there is a 
$\bP^1$
intersecting each one of the six components $|M|_a^+$, $|M|_a^-$,
$|M|_b^+$, $|M|_b^-$, $|M|_{a+b}^+$, $|M|_{a+b}^-$, in a point. See section
\ref{Geoinvest}. 
\v8

Our next result describes $Z_k(\S)$ (resp. $Z'_k(\S)$) as
representations of the group $\E$ (resp. $\E'$). (Here, $\E$ acts on
$Z_k(\S)$ {\em via} the natural action
of $\rho_a^{\otimes k}$ on $\L^{\otimes k}$, and similarly for
$\E'$.) This is based on the following result  obtained by
applying the Lefschetz-Riemann-Roch formula. We assume $\alpha\in
J^{(2)}$ is non-zero, and $2a=\a$. 

\begin{theorem}\label{1.2} One has 
$$ \tr(\rho_a^{\otimes k}) = 
\frac{1+(-1)^k}{2}\left(\frac{k+2}{2}\right)^{g-1}$$
and (for even $k$)
$$\tr({\rho'_\alpha}^{\otimes k/2}) = (-1)^{k/2}\left(
  \frac{k+2}{2}\right)^{g-1}$$
\end{theorem}

\noi{\bf Note.} In the twisted case and for
levels divisible by $4$, this result is due to Pantev \cite{Pa}. In
the untwisted case and for
levels divisible by $4$, it is also contained in
Beauville's recent paper 
\cite{Be3}. Our computation was done independently of his. 
\v8

The characters of  $Z_k(\Sigma)$ and $Z'_k(\S)$ (as representations of
the groups $\E$ and $\E'$, respectively) are determined  by the
formulas in theorem \ref{1.2} together with the trace of the
identity, given by the Verlinde formulas. Therefore the above result 
 determines  these representations up to isomorphism.

A remarkable consequence of  this is the following theorem, which was 
actually the main motivation for the present
paper. 

\begin{theorem}\label{2.4} The spaces $Z_k(\Sigma)$ and $Z'_k(\S)$
are isomorphic, as representations of the groups $\E$ and $\E'$, 
respectively, to the
TQFT-vector spaces constructed in
\cite{BHMV}.
 \end{theorem} 

This  will be discussed in more detail  in section
\ref{taugamma}.

\begin{corollary} \label{2.5} If the level $k$ is even, the spaces $Z_k(\Sigma)$
  and $Z'_k(\S)$ are decomposed as direct sums of isotypic components 
  (called `bricks' in what follows) 
  for the action of $\E$ and $\E'$. If $k\equiv 0$ mod $4$, the bricks
  are  indexed by characters of  $J^{(2)}$. If $k\equiv 2$ mod $4$, the bricks
  are  indexed by $\theta$-characteristics on the curve $\S$ (or,
  equivalently, by spin structures on $\S$). If the level is odd,
  $Z_k(\Sigma)$ is isomorphic, as representation of the group $\E$, to
  a direct sum of copies of $Z_1(\S)$ or to a direct sum of copies of the
  conjugate representation, $\overline{Z_1(\S)}$, according to
  the parity of $(k-1)/2$. 
\end{corollary} 

\begin{remarks}{\em (i) Note that formula (\ref{rhoa}) implies 
\begin{equation} \rho_a^{\otimes k} \,\rho_b^{\otimes k}\,=
\,\lambda_4(a,b)^k\,\rho_{a+b}^{\otimes k}. \label{rhoak} 
\end{equation} Hence the action
of $\E$ on
$Z_k(\S)$ factors through an action of $\E'$ if $k$ is even, and
furthermore through
an action of $J^{(2)}$ if $k\equiv 0$ mod $4$. This explains why the
bricks are indexed by characters of  $J^{(2)}$ if $k\equiv 0$ mod
$4$. If $k\equiv 2$ mod $4$, the  index set are  the characters of
$\E'$ which do not
factor through $J^{(2)}$; such characters can be identified with
$\theta$-characteristics.

(ii) The dimensions of the bricks are the same for all non-trivial
characters of 
$J^{(2)}$ in the case $k\equiv 0$ mod $4$, and depend only on the
parity of the
$\theta$-characteristic in the case $k\equiv 2$ mod $4$. In  section
\ref{3}, we will give explicit formulas for their dimensions.}
\end{remarks}

\noi{\bf Note.} Our brick decomposition can be viewed as a
generalisation of an old result of Beauville's \cite{Be2} (see also
Laszlo \cite{L}) at level $2$. Beauville constructed bases of $Z_2(\S)$ (resp. 
$Z'_2(\S)$)  whose basis elements are indexed by even (resp. odd)
$\theta$-characteristics.  This corresponds from our point of view to
the fact that  the bricks in level
$2$ are zero, if the $\theta$-characteristic has the `wrong'
parity, and one-dimensional otherwise. Note, however, that the situation in
level $2$ is very special; in general, the bricks are non-zero for
both parities. We would like to mention also
the work of van Geemen and Previato \cite{vGP1,vGP2}, Oxbury and Pauly
\cite{OP}, and Pauly \cite{P}, whose work  contains a description of
the bricks in  level $4$.

\begin{remark}\label{logstruct}{\em The proofs of our results are
    organized as follows. After a description of the fixed
  point varieties in section \ref{fpv}, theorem
  \ref{evencasei} concerning their triple intersections  is proved in section \ref{triple}. As already
  observed, theorem \ref{evencasei} implies part of theorem \ref{1.1}; the remainder of the
  proof of \ref{1.1} is given in section \ref{Geoinvest}, using some
  results about the Hecke correspondence proved in section \ref{Heckecorr}. 
  Finally, theorem 
  \ref{1.2} is proved  in section \ref{tracecomp}.  Both section
  \ref{taugamma} (where theorem \ref{2.4} is explained)  and section
  \ref{3}  (where  the brick decompositions in 
  corollary \ref{2.5} are derived) 
  require only the statements of theorem \ref{1.1} and theorem
  \ref{1.2}, but not their proofs.
}\end{remark}

\section{Relationship with the \cite{BHMV}-theory.}\label{taugamma}

In
\cite{BHMV}, the Kauffman bracket (at a
$2p$-th root of unity called $A$) was used to give a combinatorial-topological
construction of TQFT-functors
$V_p$ on a certain $2+1$-dimensional cobordism category.  It is
expected that 
these functors for $p=2k+4$
correspond (in some natural sense) to
Witten's ones for $SU(2)$ and level $k$. For example, the
dimensions of the vector space $V_p(\Sigma)\otimes \C$ \footnote{By 
definition, $V_p(\Sigma)$ is a module over a
certain abstract cyclotomic ring $k_{p}$. By $V_p(\S)\otimes \C$, we mean the
vector space obtained by extending coefficients from $k_p$ to $\C$;
this depends on a 
choice of a $2p$-th root of unity $A$ in $\C$.} associated to a closed
oriented 
surface $\S$ is also  given by the Verlinde
formula (\cite{BHMV}, cor. 1.16).

The aim of the present section is to explain the following more
precise statement of theorem \ref{2.4}. (As already mentioned in the
introduction, the reader interested only in
the algebraic geometry may proceed directly to section \ref{3}.)

\begin{theorem}\label{2.1} Let $p=2k+4$. For the `right' choice of
  $2p$-th root of unity $A$, the vector spaces $V_p(\S)\otimes \C$ 
 and
  $Z_k(\S)$ are isomorphic as representations of $\E$
  (this group is denoted by $\Gamma(\S)$ in \cite{BHMV}). This
 isomorphism sends the involution $\rho_a^{\otimes k}$ acting on
 $Z_k(\S)$ to the involution $\tau_\gamma$ acting on $V_p(\S)\otimes
 \C$, where $\g$ is a simple closed curve on $\S$ representing the
 Poincar\'e dual of $a$.  A similar statement holds in the twisted
case.
\end{theorem}

\noi{\bf Note.} 
For
simplicity of exposition, a few technical details related to framing
issues will be
suppressed from the discussion here. The reader is also referred to the
survey article \cite{MV}.
\v8

By definition, elements of  $V_p(\Sigma)$
are represented by linear combinations of compact oriented
$3$-manifolds 
$M$ with boundary
equal to $\S$. (No complex structure on $\S$ is needed here.) The
$3$-manifolds may also contain colored links or, more generally,
colored trivalent graphs. In the following, we
assume $p$ is even and put  $p=2k+4$. Then a color is
just an element of $\{0,1,2,\ldots,k\}$. 

In the \cite{BHMV}-theory, one has involutions $\tau_\gamma$ of
$V_p(\Sigma)$ associated to unoriented simple closed curves $\gamma$ on
$\S$. They are defined as follows. Consider a
vector represented by some $(M,L)$, where $M$ is a $3$-manifold with
boundary $\S$, and $L$ stands for 
some colored link inside $M$. Then the action of  $\tau_\gamma$
consists of adding to the link $L$ already present in $M$, an
additional component consisting of the curve
$\gamma$ pushed slightly inside $M$ in a neighborhood of $\S=\partial
M$, where $\gamma$ is  colored by $k$ (the last
color). It was shown in
section 7 of \cite{BHMV} that (because of the relations which hold in
$V_p(\S)$) this endomorphism $\tau_\gamma$ is an
involution.

Following  \cite{BHMV}, p. 917, the groups generated by these
involutions 
acting on
$V_p(\Sigma)$ can be described as follows. Let $\Gamma(\S)$ be the
group with one generator $[a]$ for each $a\in H_1(\S;\Z/4)$ plus one
additional generator $u$, and the following relations: the element $u$
is central, $u^4=1$, $[a]^2=1$ for all $a\in H_1(\S;\Z/4)$, and 
\begin{equation}
[a]\,[b]\, = \,
u^{a\cdot b} \,[a+b]\label{uu}
\end{equation}
 for all $a,b\in H_1(\S;\Z/4)$. Here, $a\cdot
b\in \Z/4$
  denotes the {\em mod} $4$ intersection form determined by the orientation
  of $\S$.

If a {\em mod} $2$ class $\alpha\in H_1(\S;\Z/2)$ is given, every lift
of $\alpha$
to a {\em mod} $4$ class $a\in H_1(\S;\Z/4)$ determines an element $[a]\in
\Gamma(\S)$. However, up to multiplication by $u^2$, this element
$[a]$ depends only  on
$\alpha$. (Indeed, if $a_1,a_2$ are two such lifts, then $[a_1]=u^{2
  a_2\cdot x}[a_2]$
where $a_1-a_2=2x$.) From this it follows
easily that the group
$\Gamma(\S)$ is a central extension of $H_1(\S;\Z/2)$ by
$\Z/4$. (It can  be viewed as some kind of `reduced' Heisenberg group
associated to twice the intersection form.)

Here is the relationship of the group $\Gamma(\S)$ with the involutions
$\tau_\gamma$. Given $\alpha \in H_1(\S;\Z/2)$, a lift of $\alpha$ to
an element of
$\Gamma(\S)$ can also be specified by an unoriented simple closed
curve $\gamma$ representing $\alpha$, as follows. The curve   
$\gamma$ determines an
element $a \in H_1(\S;\Z/4)$ up to sign, and the induced element
$[a]\in \Gamma(\S)$ is well-defined,
since $a\cdot a=0$ and hence $[-a]=[a]$. 

It is clear from the geometric description of
the $\tau_\gamma$'s given above  that the commutation properties of
these involutions are related to the
intersection properties  of the associated curves; for example two
such 
involutions 
obviously commute if the corresponding simple closed curves on $\S$
don't intersect. The fact that the 
$\tau_\gamma$'s satisfy precisely the relations in $\Gamma(\S)$ is
shown in prop. 7.5 
of \cite{BHMV}. Moreover, it is shown there  that  the
central
 element $u$ acts on $V_{2k+4}(\S)$ as multiplication by
 $(-1)^{k+1}A^{(k+2)^2}$.

Now let us identify  $\Z/4=\mu_4$ such that $1\in \Z/4$ corresponds to
$i\in\mu_4$. Also, let us use Poincar\'e duality to identify $
J^{(4)}(\S)=H^1(\S;\mu_4)$ with $H_1(\S;\Z/4)$. By a well-known folk
theorem, there exists a
sign $\varepsilon =\pm 1$ such that the Weil pairing  is
related to the intersection form as follows: \footnote{The value of
  $\varepsilon$ should be  known, but  we have been
  unable to locate it in the
  literature.}  $$
\lambda_4(a,b)=(\varepsilon i)^{a\cdot b}.$$ We now choose the
primitive root of unity $A$ of order $2p=4k+8$ such that 
$(-1)^{k+1}A^{(k+2)^2}=(\varepsilon i)^k$. For instance, we may choose 
$A=-\varepsilon
e^{2i\pi/(4k+8)}$.   Comparing (\ref{rhoa}) and
  (\ref{uu}), we see that the assignment
  $$\tau_\gamma \mapsto \rho_a^{\otimes k}$$ (where the curve $\gamma$
  represents the Poincar\'e dual of $a$) defines an isomorphism
  from the image of $\Gamma(\S)$ in the general linear
group of $V_{2k+4}(\S)\otimes \C$ to the image of the group $\E$ in
  the general linear
group of $Z_k(\S)$. (In particular, for $k=1$ we have an isomorphism
  $\Gamma(\S)\cong \E$.) 

Using the
dimension formulas given in section 7 of \cite{BHMV}, one can check 
that the traces of the involutions $\tau_\gamma$ acting on
$V_p(\Sigma)$ coincide precisely with the traces computed in our theorem
\ref{1.2}. This verification will be omitted. 

Since representations of finite groups are determined by their
characters, this proves theorem \ref{2.1} in the untwisted case. 

 The twisted
case,  where $Z_k(\S)$ is replaced with $Z'_k(\S)$, is similar. Here
we must replace $V_p(\S)$  with $V_p(\S')$, where $\S'$ is $\Sigma$
with one puncture
colored by $k$ (the last color). (We maintain the convention
$p=2k+4$.) Then we have again an isomorphism of representations (the
group corresponding to $\E'$ is called $\Gamma'(\S)$ in
\cite{BHMV}). The dimension of
$V_p(\S')$ is given by the `twisted Verlinde formula', see
\cite{BHMV}, Remark 5.11. This space is zero if $k$ is odd, and this
is why we define $Z'_k(\S)$ to be zero for odd $k$. (See also Thaddeus
\cite{Th}.) (N.b., the vector
spaces $V_p$ 
are also defined for surfaces with colored punctures; they are called
`surfaces with colored structure' in \cite{BHMV}. The theory of the involutions
$\tau_\gamma$ works the same for these, with the {\em caveat} that the
curves $\gamma$ must avoid the punctures. We expect analogues of our
results to hold in the case of general colored punctures, using moduli spaces
of parabolic bundles.)

\v8
\noi{\bf Remark.} Although it has been known for some time  that $V_p(\S)$ and
  $Z_k(\S)$ have the same dimensions, it seems, to the best of our
  knowledge, that a {\em natural} isomorphism between the two 
theories is still missing. Of course, the fact that these two spaces
are isomorphic also as representations (of canonically isomorphic groups) gives further evidence that there should be a
natural isomorphism  between the two theories.

\section{Brick decomposition and dimension formulas.}\label{3}

In this section, we describe $Z_k(\S)$ and $Z'_k(\S)$ as
representations of the groups $\E$ and $\E'$, respectively. We also
give explicit dimension formulas. The results of  this section are
immediate consequences of the isomorphism of theorem \ref{2.1} and the
computations in  \cite{BHMV} and \cite{BM}. In order to make this
paper self-contained, we will, however, derive them directly from
theorems \ref{1.1} and \ref{1.2}.

We first recall the Verlinde formula and its twisted counterpart. Put $d_g(k)=dim\ Z_k(\S_g)$ and $d'_g(k)=dim\
Z'_k(\S_g)$, where the subscript $g$ indicates the genus of the curve $\S_g$. \footnote{These numbers are denoted by $d_g(p)$
and $\widehat d_g(p)$ in \cite{BHMV,MV,BM}, where $p=2k+4$.} Then one has  
$$d_g(k)=\left(\frac{k+2}{2}\right)^{g-1}\
\sum_{j=1}^{k+1} \ \left(
\sin \frac{ \pi j}{k+2}\right)^{2-2g}$$

$$d'_g(k)=\left(\frac{k+2}{2}\right)^{g-1}\
\sum_{j=1}^{k+1} \ (-1)^{j+1}\left(
\sin \frac{ \pi j}{ k+2}\right)^{2-2g}$$

\noi{\bf The case $k \equiv 0$ mod $4$.}

In view of formula
(\ref{rhoak}) in section \ref{smr}, the action of the group $\E$
on $Z_k(\S_g)$ factors through an action of
$J^{(2)}(\S_g)=H^1(\S_g;\mu_2)$.  We
have a direct sum decomposition 

$$Z_k(\S_g)=\bigoplus_{h} Z_k(\S_g;h)$$
where $h$ runs through the characters of $J^{(2)}$, and $Z_k(\S_g;h)$
is the subspace of $Z_k(\S_g)$ where $\a$ acts as multiplication by
$h(\a)\in \mu_2={\pm 1}$, for all  $\alpha\in J^{(2)}$. We will refer
to $Z_k(\S_g;h)$ as the {\em brick} associated to $h$. 

 By theorem \ref{1.2}, the character of the representation $Z_k(\S_g)$
 takes the same value on all non-trivial elements of $\Jt$, and is
 therefore invariant under the automorphism group of $\Jt$. It follows
 that the dimension  of the brick $Z_k(\S_g;h)$ is the
same for all non-trivial characters $h$ of the group, since the
 automorphism group acts transitively on the set of non-trivial
 characters.  We denote this dimension  by $d_g^{(1)}(k)$,
and we put  $d_g^{(0)}(k)=dim \
Z_k(\S_g;0)$, where $0$ denotes the trivial character. (Thus, the
space  $Z_k(\S_g;0)$ is  the fixed point set of the action of the
group $J^{(2)}$ on $Z_k(\S_g)$.)  These numbers can be computed from the following two formulas:
$$d_g(k)=d_g^{(0)}(k) +(2^{2g}-1)\, d_g^{(1)}(k)$$
$$d_g^{(0)}(k)-d_g^{(1)}(k)=tr(\alpha)=
\left(\frac{k+2}{2}\right)^{g-1}, \ \ \alpha\neq 0$$ For example, one
has $$d_g^{(0)}(k)=\frac{1}{2^{2g}} \left(d_g(k) + (2^{2g}-1) \left(\frac{k+2}{2}\right)^{g-1}\right)$$

Similarly, $Z'_k(\S_g)$ is the direct sum of  bricks 
$Z'_k(\S_g;h)$, 
and the dimensions ${d'_g}^{(0)}(k)=dim \ Z'_k(\S_g;0)$ and ${d'_g}^{(1)}(k)=dim \ Z'_k(\S_g;h)$ for $h\neq 0$, can be computed as before (just replace $d_g(k)$ with $d'_g(k)$ and $d_g^{(\varepsilon)}(k)$ with ${d'_g}^{(\varepsilon)}(k)$ in the above).

\v8\noi{\bf Example: The case $k=4$.} The numbers $d_g^{(1)}(4)$ and
${d'_g}^{(1)}(4)$ are equal to $(3^{g-1}+1)/2$ and $(3^{g-1}-1)/2$,
respectively, and the numbers $d_g^{(0)}(4)$ and ${d'_g}^{(0)}(4)$
are obtained by adding $3^{g-1}$. These numbers have appeared in Oxbury and Pauly
\cite{OP} and Pauly \cite{P}.
\v8
\noi{\bf The case $k\equiv 2$ mod $4$.}\nopagebreak 

In this case, the action of
$\E$ factors through an action of the group $\E'$ but not through an
action of $J^{(2)}$. Indeed, the involutions $\rho_a^{\otimes 2}$ acting on $\L^2$ depend only on
$\alpha$, and satisfy the same relations as the
$\rho'_\alpha$'s (see formulas (\ref{rhoprima}) and 
(\ref{rhoak}) in section \ref{smr}, and note that $\lambda_4(a,b)^2=\lambda_2(\a,\b)$). Therefore one has a direct sum
decomposition $$Z_k(\S_g)=\bigoplus_{q} Z_k(\S_g;q)$$ where $q$ runs
through the characters of the group $\E'$ which do not factor through
$J^{(2)}$, {\em i.e.}, such that $q$ takes the value $-1$ on the central
element $-1\in\mu_2\subset \E'$. Such characters are in $1$-to-$1$ correspondence with
functions $q\colon J^{(2)}\rightarrow \mu_2$  such that
$$q(\alpha+\beta)\,=\,q(\alpha)\,q(\beta)\, \lambda_2(\alpha,\beta).$$ In
other words, $q$ runs through the set of quadratic forms on
$J^{(2)}\cong H^1(\S_g;\Z/2)$
inducing the  Weil pairing, {\em i.e.}, the {\em mod} $2$ intersection form. It
is well-known  that such quadratic forms correspond to spin
structures, or equivalently, $\theta$-characteristics,  on $\S$. (See
Atiyah \cite{At2}, Johnson \cite{Jo}.)

 Let ${\G}\cong Sp(2g;\Z/2)$ be the group of automorphisms of $\Jt$ preserving the
order $2$ 
Weil pairing $\lambda_2$. The group $\G$ acts on $\E'$. As in the case
$k\equiv 0$ mod $4$, it follows
from theorem \ref{1.2} that the character of
  the representation $Z_k(\S_g)$ is invariant under the
  action of $\G$. It is well-known that the induced action of $\G$ on
  quadratic forms  has two orbits which are characterized by the {\em Arf
  invariant}, {\em i.e.}, the action is transitive on the set of forms $q$ with the
same Arf invariant $\Arf(q)\in \Z/2$. (The Arf invariant of the quadratic form
corresponds to the parity of the $\theta$-characteristic.) This shows
that  the
  dimension of the brick $Z_k(\S_g;q)$ depends only on $\Arf(q)$. Put $d_g^{(\varepsilon)}(k)=dim \ Z_k(\S_g;q_\varepsilon)$, where $q_\varepsilon$ has Arf invariant $\varepsilon\in \Z/2$. These dimensions can be computed  from the following two formulas:
$$d_g(k)=2^{g-1}(2^g+1)\,d_g^{(0)}(k) +2^{g-1}(2^g-1)\, d_g^{(1)}(k)$$
$$2^{g-1}(d_g^{(0)}(k)-d_g^{(1)}(k))=tr(\rho_a^{\otimes k})=
\left(\frac{k+2}{2}\right)^{g-1}, \ \ \alpha\neq 0$$

The first formula follows from the fact that the number of quadratic
forms with zero Arf invariant is equal to $2^{g-1}(2^g+1)$. The second
formula follows from the fact\footnote{A nice way to think about this
  is to observe that there is a natural bijection between quadratic
  forms $q$ with fixed Arf invariant and such that $q(\a)=-1$, and
  quadratic forms with arbitrary Arf invariant on the
  $2g-2$-dimensional space $\lag \a\rag ^\bot/\lag \a\rag$.} that for $\alpha\neq 0$, one has 
$$\sharp \{q\,|\,q(\alpha)=-1,\ \Arf(q)=0\}=2^{2g-2}=\sharp \{q\,|\, q(\alpha)=-1,\ \Arf(q)=1\}$$

Similarly, $Z'_k(\S_g)$ is the direct sum of  bricks  $Z'_k(\S_g;q)$,
and the dimensions ${d'_g}^{(\varepsilon)}(k)=dim \
Z'_k(\S_g;q_\varepsilon)$ can be computed as before (just replace
$d_g(k)$ with $d'_g(k)$ and $d_g^{(\varepsilon)}(k)$ with
${d'_g}^{(\varepsilon)}(k)$ in the above, but notice that
$tr({\rho'_\a}^{\otimes k/2})$ is now equal to $-\left((k+2)/2\right)^{g-1}$.)

Here are explicit formulas for the dimensions. They are equivalent to  the
formulas on p. 264 of \cite{BM} \footnote{Warning: Note that $k$ does not
denote the level in \cite{BM}; one has $p=8k$ in \cite{BM} while in
the present paper $p=2k+4$.}.
\begin{eqnarray} d_g^{(\varepsilon)}(k)&=&{\frac{1}{2^{2g}}} \left(d_g(k)+\left(\frac{k+2}{2}\right)^{g-1} \left((-1)^{\varepsilon} 2^g-1\right)\right)\nonumber\\
{d'_g}^{(\varepsilon)}(k)&=& {\frac{1}{2^{2g}}}\left(d'_g(k) +\left(\frac{k+2}{2}\right)^{g-1} \left(1-(-1)^{\varepsilon} 2^g\right)\right)\nonumber
\end{eqnarray}

\v8\noi{\bf Example: The case $k=2$.} The numbers $d_g^{(0)}(2)$ and
${d'_g}^{(1)}(2)$ are equal to $1$,  and the numbers $d_g^{(1)}(2)$
and ${d'_g}^{(0)}(2)$ are zero. Therefore  $d_g(2)=2^{g-1}(2^g+1)$ and
$d'_g(2)=2^{g-1}(2^g-1)$ (see Beauville \cite{Be2}).

\v8\noi{\bf Note.} In \cite{BM}, a $\Z/2$-graded TQFT-functor is
constructed on a cobordism category of surfaces equipped with spin
structures (and other things).  Given a connected surface $\S$ with a
spin structure $\sigma$, the even (resp. odd) part of this functor is
isomorphic to  $Z_k(\S_g;q_\sigma)$ (resp. $Z'_k(\S_g;q_\sigma)$),
where $q_\sigma$ is the quadratic form corresponding to $\sigma$.

\begin{remark}{\em The action of the symplectic group $\G\cong Sp(2g;\Z/2)$
    permuting the bricks has the following geometric interpretation.  Recall that all elements of $\G$ can be represented by 
    diffeomorphisms of $\S$. In the \cite{BHMV}-theory, one has, more or less by
    definition, a
    natural 
action of a
    certain extended mapping class group  on $V_{2k+4}(\S)$. On the
    geometric side, there is also  a (projective-linear) action
    of the mapping class group of $\S$ on $Z_k(\S)$; this action is constructed using
    Hitchin's projectively-flat
    connection
 \cite{Hi}. It is, of course, expected that $V_{2k+4}(\S)\otimes \C$ and $Z_k(\S)$ are isomorphic as
representations of the extended mapping class group. In any case, it
    is easy to see in both theories that  the
action of 
a diffeomorphism
$f$  takes the brick associated to a character $h$ (resp. a quadratic
form $q$)  to the brick associated to $f^*(h)$ (resp.
$f^*(q)$). On the geometric side, the main reason for this is  that  Hitchin's 
    connection is (projectively) invariant under the actions of both
    the mapping class group and the group $\E$.
}\end{remark}

\noi{\bf The odd-level case.} 

\begin{proposition} If $k\equiv 1$ mod $2$,
  $Z_k(\Sigma)$ is isomorphic, as representation of the group $\E$, to
  a direct sum of copies of $Z_1(\S)$, if $k\equiv 1$ mod $4$, and to a direct sum of
  copies of the
  conjugate representation, $\overline{Z_1(\S)}$, if $k\equiv 3$ mod
  $4$.  
\end{proposition}
\proof  Note that the character
of the representation $Z_1(\S)$ is the function $\chi\colon \E\ra \C$ 
which is zero on all group elements not in the central subgroup $\mu_4\subset
\E$, while  the trace of a central element $\lambda\in \mu_4$
is $\chi(\lambda)=2^g\,\lambda$. The character of the conjugate
representation, $\overline{Z_1(\S)}$, is of course the conjugate
character $\overline{\chi}$. Now it follows from theorems \ref{1.1}
and \ref{1.2} and formula (\ref{rhoak}) that
the character of $Z_k(\S)$ is a multiple of $\chi$ if $k\equiv 1$ mod $4$, and a multiple of $\overline{\chi}$, if $k\equiv 3$ mod
  $4$. It is however well-known that $\chi$ is an irreducible
  character, and this proves the proposition. 

\v8
\noi{\bf Note.} This result corresponds, {\em via} the isomorphism
$Z_k(\S)\cong V_{2k+4}(\Sigma) \otimes \C$, to theorems 1.5
  and 1.6 of
  \cite{BHMV}. It is shown there that (for odd $k$) 
$$V_{2k+4}(\Sigma)\cong V_2'(\Sigma)\otimes
V_{k+2}(\Sigma)$$ as representations of the group $\E$,  where
$V_2'(\Sigma)$ and $V_{k+2}(\Sigma)$ are defined in
\cite{BHMV}. Moreover, 
the group $\E$ acts trivially on $V_{k+2}(\Sigma)$, and
$V_2'(\Sigma_g)$ is isomorphic to $V_6(\Sigma)$, and hence to $Z_1(\S)$
or to $\overline{Z_1(\S)}$, after a change of
coefficients. It would be interesting to have an algebro-geometric
interpretation of the space $V_{k+2}(\Sigma)$, which in the
\cite{BHMV}-theory can be
interpreted as a
$SO(3)$-TQFT vector space. Here, the name $SO(3)$-TQFT just means that the allowed colors in
this TQFT are even, or in other words, are $SU(2)$-representations
which lift to $SO(3)$.

\section{The fixed point varieties.}\label{fpv}

We need to analyze the action of
$\Jt$ on $M$ and $M'$ in order to obtain information about our
lifts. In this section, we describe  the various fixed
point varieties, mainly following Narasimhan and Ramanan
\cite{NR}. We also discuss 
when two lifts of $\a\in \Jt$ to $\Jf$ determine the same component of
the fixed point set $|M|_\a$. 

 \v8\noi{\bf Notation.} Throughout this paper, we denote by
 $L_\a$ (resp. $L_a$) the line bundle on
$\S$ corresponding to $\a\in
\Jt$ (resp. $a\in \Jf$). 
\v8

Let $\a\in\Jt$ be nonzero. Let $\pi_\a : \Sigma^\alpha \rightarrow \Sigma$ be the 2-sheeted unramified covering of
$\Sigma$ corresponding to
$\alpha$, and let $\phi_\a$ be the covering transformation  of $\pi_\a$.
Using the line bundle $L_\a$, we can explicitly construct $\S^\a$ as
$$\S^\a = \left\{ \xi \in L_\a | \xi\otimes\xi = 1\right\}$$
using an isomorphism $L^2_\a\cong \cO_\S$. The involution $\phi_\a$ is then
induced by multiplication by $-1$ on $L_\a$.

Given a line bundle $L$ over $\S^\a$, the push-down $E=\pia(L)$ can be
obtained by descending $L\oplus
\phi_\a^*(L)$ (which is naturally an equivariant bundle) to $\S$:
\begin{equation}
E=(L\oplus \phi_\a^*(L)){\big/}\langle \phi_\a\rangle.\nonumber
\end{equation} Here, the
natural involution of $L\oplus \phi_\a^*(L)$ covering $\phi_\a$ is
again denoted by $\phi_\a$.  The
fundamental observation is that $E\otimes L_\a$ is 
isomorphic to $E$.  This follows formally from the pull-push formula
$\pia(L)\otimes L'\cong \pia(L\otimes \pi_\a^*L') $ and the fact
that $L_\a$ pulls back to the trivial bundle on $\S^\a$. 

\begin{remark} \label{fundobs} {\em We will later need the following
  explicit isomorphism from $E=\pia(L)$ to $E\otimes L_\a$. Note that $\pi_\a^*L_\a$, as an equivariant bundle, is
isomorphic to $\cO_{\S^\a}^-$, that is, the  trivial line bundle $\S^\a\times
\C$, but with non-trivial action, given by $(x,z)\mapsto (\phi_\a(x),
-z)$. Therefore $$E\otimes L_\a=(L\oplus
\phi_\a^*(L)){\big/}\langle -\phi_\a\rangle.$$ It follows that the diagonal automorphism
$1\oplus(-1)$ of $L\oplus \phi_\a^*(L)$ descends to an isomorphism from $E$
to $E\otimes L_\a$.}\end{remark}

Let $\Nm_\a\colon \Pic(\Sigma^\alpha) \rightarrow 
\Pic(\Sigma)$ be the  classical albanese or norm map 
induced by the following map on divisors. If $D= \sum n_j x_j\in \Div(\S^\a)$ then
$\Nm_\a(D) = \sum n_j \pi_\a(x_j)\in \Div(\S)$.

The {\em Prym variety} $P_\alpha$ associated to $\alpha$ is by definition the
connected component of $\Nm_\a^{-1}(\cO_\S)$ containing
$\cO_{\Sigma^\alpha}$. It is a principally polarized abelian variety
of dimension $g-1$. The quotient variety
$P_\alpha/\langle\pm 1\rangle$ is  called the {\em Kummer variety.} 

 We define $\theta_\a : \Pic(\Sigma^\alpha) \rightarrow
\Pic(\Sigma)$ by $$\theta_{\a}(L) = \det(\pia(L)) =\Nm_\a(L)\otimes
L_\alpha$$
(see \cite{NR} for the second equality).  Note that  
$\phi_\a$ acts on 
$\Pic(\Sigma^\alpha)$ by sending $L$ to $\phi_\a^*(L)$.

\begin{proposition}[Narasimhan and Ramanan \cite{NR}]\label{6.1} (i) The map $L\mapsto (\pi_\a)_*(L)$
  induces  isomorphisms 
\begin{equation}\label{thetu}
\theta_\a^{-1}(\cO_\S)/\lag\phi_\a \rag \,\mapright\sim \,|M|_\a
\end{equation} 
\begin{equation}\label{thett}\theta_\a^{-1}([p])/\lag\phi_\a\rag \,\mapright\sim \,|M'|_\a
\end{equation}
(ii) Moreover, $\theta_\a^{-1}(\cO_\S)/\lag\phi_\a\rag$ is isomorphic to two
copies of $P_\alpha/\lag\pm 1\rag$, while $\theta_\a^{-1}([p])/\lag\phi_\a\rag$ is
isomorphic to $P_\a$.
\end{proposition}

\begin{remark}{\em In the twisted case, 
    $\pia\colon \theta_\a^{-1}([p]) \ra |M'|_\a$ is a double covering,
    with covering transformation $\phi_\a$. In the untwisted case, the
    same holds on the open subvariety of $\theta_\a^{-1}(\cO_\S)$ 
where $\phi_\a$ acts freely, while the fixed
points of $\phi_\a$ are sent by the map $\pia$ bijectively to the points in $|M|_\alpha$
represented by semi-stable but not stable bundles. (See remark
\ref{check} below.)
}\end{remark}

 For later use, we need the following more explicit
description of $\theta_\a^{-1}(\cO_\S)$ and $\theta_\a^{-1}([p])$. For $d=0,1$, define $ \Phi^d_\a :  \Pic_d(\S^\a) \ra \Pic_0(\S^\a)=J(\S^\a)$
by 
\begin{equation} \Phi^d_\a(L) = L\otimes\phi^*_\a L^{-1}.\nonumber
\end{equation} The following two properties (\ref{imPhi1}) and
(\ref{imPhi2}) of the maps $
\Phi^d_\a$ are elementary facts from the classical theory of
line bundles on curves, see {\em e.g.} Appendix  B in \cite{ACGH}. 

First,  one has the disjoint union 
\begin{equation} \Nm_\a^{-1}(\cO_\S)=\im\Phi^0_\a \cup \im \Phi^1_\a
    \label{imPhi1}
\end{equation} and the Prym
variety $P_\a$ is  equal to the  component $\im\Phi^0_\a$. 

Second, note that $\Nm_\a(\pi_\a^*(L_\b))=L_\b^{\otimes 2}=\cO_\S$ for all
$\beta\in \Jt$. Moreover, one has 
\begin{equation} 
\pi_\a^*L_\b \in \im \Phi^d_\a \ \Leftrightarrow\
\lambda_2(\a,\b)=(-1)^d. \label{imPhi2}
\end{equation}

Here, $\lambda_2 : J^{(2)}\times J^{(2)} \ra \mu_2$ is  the order $2$ Weil
pairing. 

Now pick  $a\in \Jf$ such that
  $2a=\a$, and $\b\in\Jt$ such that $\lambda_2(\a,\b)=-1$.  Note that
  $a'=a+\b$ is another element of $ \Jf$ such that $2a'=\a$. Also, pick
  a point 
  $\pa\in\pi_\a^{-1}(p)\subset \S^\a$.  Note that
$\theta_\a^{-1}(\cO_\S)$ and $\theta_\a^{-1}([p])$ are both isomorphic
to $\Nm_\a^{-1}(\cO_\S)$. {}From (\ref{imPhi1}) and
  (\ref{imPhi2}), we have the following description of
  $\theta_\a^{-1}(\cO_\S)$ and $\theta_\a^{-1}([p])$ as disjoint
  unions:
\begin{eqnarray} \theta_\a^{-1}(\cO_\S)&=&\pi_\a^*L_{a}\otimes
  \im\Phi^0_\a \,\cup\, \pi_\a^*L_{a'} \otimes \im \Phi^0_\a \label{untw1}\\
\theta_\a^{-1}(\cO_\S)&=&\pi_\a^*L_{a}\otimes
  \im\Phi^0_\a \,\cup\, \pi_\a^*L_{a} \otimes \im \Phi^1_\a \label{untw2}\\
\theta_\a^{-1}([p])&=&\pi_\a^*L_{a}\otimes[p_\a]\otimes
  \im\Phi^0_\a \,\cup\, \pi_\a^*L_{a} \otimes [\phi_\a(p_\a)]\otimes\im
  \Phi^0_\a \label{tw}
\end{eqnarray}

\begin{remark} {\em  It follows from (\ref{untw2}) and (\ref{tw}) that the action of 
$\phi_\a$ preserves the two components of $\theta_\a^{-1}(\cO_\S)$,  while it
exchanges the two components of $\theta_\a^{-1}([p])$. By (\ref{imPhi1}), the action of 
$\phi_\a$ on
$\Pic(\Sigma^\alpha)$  restricts to
multiplication by $-1$ (that is, the map  $L\mapsto L^{-1}$) on
$\Nm_\a^{-1}(\cO_\S)$. This proves \ref{6.1}(ii).
}\end{remark} 

We can now describe  
when two lifts of $\a\in \Jt$ to $\Jf$ determine the same component of
the fixed point set $|M|_\a$. Recall that $|M|_a^+$ was defined to be 
the component of $|M|_\a$ containing the S-equivalence
class of the semi-stable bundle $L_a\oplus L^{-1}_a$.

\begin{proposition}\label{4.6} (i) Let $a_1,a_2 \in \Jf$ such that
  $2a_1=2a_2=\a$. Then $$|M|_{a_1}^+
= |M|_{a_2}^+ \ \Leftrightarrow\ \lambda_2(a_1-a_2,\a)=1\in\mu_2.$$

(ii) For $\b\in \Jt$, the action of 
$\b$ on $|M|_\a$ interchanges the two components of $|M|_\a$ if
and only if  $\lambda_2(\a,\b)=-1$.
\end{proposition}

\proof Note that  $\pia(\pi_\a^*
(L_{a_i}))\cong L_{a_i}\oplus L_{a_i}\otimes L_\a\cong L_{a_i}\oplus L_{a_i}^{-1}$. Therefore $ |M|^+_{a_i}$ is the
component 
$\pi_{\a*}(\pi_\a^*L_{{a_i}}\otimes \im\Phi^0_\a)$, and (i) follows
from formula (\ref{untw1}). Now part (ii) follows  from (\ref{imPhi2}), since
the action of $\b$ on $|M|_\a$ lifts to tensoring with $\pi_\a^*L_\b$ on $ \theta_\a^{-1}
(\cO_\S)$. This completes the proof.

\v8\noi{\bf Note.} Translating the `multiplicative' notation of the
Weil pairing into the `additive' notation in  section
\ref{taugamma}, the condition $\lambda_2(a_1-a_2,\a)=1\in\mu_2$
becomes the condition $((a_1-a_2)/2)\cdot \a=0\in \Z/2$. Thus, 
prop. \ref{4.6}(i) implies that $|M|_{a_1}^+
= |M|_{a_2}^+$ if and only if $[a_1]=[a_2]$, where $[a_i]$ is the lift
of $\a$ to the group $\Gamma(\S)$ defined in section
\ref{taugamma}.

\begin{remark}\label{check}{\em  Let $M^{sing}$ denote the set of points of $M$
    represented by semistable, but not stable, bundles. A semi-stable
    bundle $E$ with $\Gr(E) \cong L_1\oplus L_2$ represents a point in
    $M^{sing}$ if and only if $L_2\cong L_1^{-1}$, and this point lies
    in $|M|_\a$ if and only if $L_1\otimes L_\a\cong L_2$. This shows
    that $M^{sing}\cap |M|_\a$ is precisely the set of points represented
    by bundles of the form $L_a\oplus L_a^{-1}$ with $a\in \Jf$ and
    $2a=\a$. 
Since $L_a\oplus L_a^{-1}\cong\pia(\pi_\a^*(L_a))$, and  $\pi_\a^*(\{L\in \Pic_0(\S) | L^2 \cong L_\a\})$ is
precisely the fixed point set of $\phi_\a$ on
$\theta_\a^{-1}(\cO_\S)$, we see that  $\pi_{\a*}$ induces a bijection from
that fixed point set to $M^{sing}\cap |M|_\a$, as asserted above.
}\end{remark}

\section{Intersections and triple intersections.} \label{triple}
\v8
In this section, we describe how the various fixed point varieties
intersect. The triple intersections  will be used in the proof of theorem \ref{1.1}.

\v8\noi{\bf Note.} The 
intersection properties of the fixed point sets in relation to the order $2$
Weil pairing are well-known; they are used for example in van Geemen and
Previato \cite{vGP1}. On the other hand, the
relationship of the triple intersections of their individual {\em components} (in the
untwisted case) with the order $4$ Weil pairing seems to be new. 
\v8

\begin{proposition}\label{6.5}
For $\alpha,\beta\in \Jt$ non-zero distinct elements, we have that
$$|M|_\alpha\cap |M|_\beta \neq
\emptyset\ \Leftrightarrow\ \lambda_2(\alpha,\beta) = 1\in\mu_2$$
$$|M'|_\alpha\cap |M'|_\beta \neq
\emptyset\ \Leftrightarrow\ \lambda_2(\alpha,\beta) = -1\in\mu_2$$
\end{proposition}

\proof
The morphism $\pi_{\a*}\colon \Pic_d(\Sigma^\alpha)\ra |M_d|_\alpha$
is surjective and $\Jt$-equivariant. (Here, $\b\in J^{(2)}$ acts on
$\Pic_d(\Sigma^\alpha)$ by $L\mapsto L\otimes \pi_\a^*L_\b$.) Hence we get the following
description of the intersection
$$|M_d|_\alpha \cap |M_d|_\beta = \pia( \{L\in \Pic_d(\S^\alpha) | L\otimes
\pi_\a^*L_\beta \cong \phi_\a^* L\mbox{ or } L\}).$$
Hence we see that
$$|M_d|_\alpha \cap |M_d|_\beta \neq \emptyset$$
if and only if 
$$\pi_\a^*L_\b \in \im \Phi^d_\a.$$
{}From  (\ref{imPhi2}), this is the case if and only if
$\lambda_2(\a,\b)=(-1)^d $. Using the action of $J(\S)$  on
$|M_d|_\alpha \cap |M_d|_\beta$, 
the results for $|M|_\alpha\cap |M|_\beta$ and $|M'|_\alpha\cap
|M'|_\beta$ follow from this.

\begin{proposition}\label{trans} If $\lambda_2(\a,\b)=1$, the quotient group $J^{(2)}/\lag\a,\b\rag$
  acts simply transitively on
  $|M|_\a\cap |M|_\b$. In particular, this 
  intersection has $2^{2g-2}$
  elements. If $\lambda_2(\a,\b)=-1$, the same holds for $|M'|_\a\cap |M'|_\b$. 
\end{proposition}

\proof Put  $I_{\a,\b}=\{L\in \theta_\a^{-1}(\cO_\S)\,|\,
\pi_\a^*(L_\b)\cong L\otimes \phi_\a
^*L^{-1}\}$. Then $\pia\colon I_{\a,\b}\ra |M|_\a\cap |M|_\b$ is a
double covering. Note that $J^{(2)}/\lag \a\rag$ acts simply 
transitively 
on $I_{\a,\b}$. The 
involution $\phi_\a$ is on $ I_{\a,\b}$ the same as  tensoring with 
$\pi_\a^*L_\b$, 
in other words, the action of $\b$. This proves the result in the
untwisted case. The twisted case is proved similarly. 
\v8
\noi{\bf Note.} In view of remark \ref{check}, this description
 shows that $|M|_\a\cap |M|_\b$ is contained in the stable part of
 $M$.
\v8
\begin{proposition}\label{count} $|M|_\a\cap |M|_\b$ is the disjoint union of the sets
 $|M|_a^\varepsilon\cap |M|_b^\mu $, 
where $\varepsilon=\pm$ and $\mu=\pm$,  each of
 which sets has
 $2^{2g-4}$ elements.
\end{proposition} 
\proof Recall from proposition \ref{4.6}(ii) that the action of $\gamma\in
\Jt$ exchanges the components of $|M|_\a$ if and only if
$\lambda_2(\a,\g)=-1$.  Thus, the result follows by exploiting  the fact
that for every choice of signs $\varepsilon=\pm 1$ and $\mu=\pm 1$, there exists $\g$
such that $\lambda_2(\a,\g) = \varepsilon$ and
$\lambda_2(\b,\g)= \mu$.
\v8
We now turn to the triple intersections. The first observation is the
following easy lemma. 
\begin{lemma} Let $\a,\b,\g$ be
distinct non-zero elements of $\Jt$ such that  the triple
intersection $|M|_\alpha \cap |M|_\beta
\cap |M|_\g$ is non-empty. Then $\g=\a+\b$. The same holds for the triple intersections
in the twisted case. 
\end{lemma} 
\proof Indeed, it
follows from the description in prop. \ref{trans} that the triple
intersection can only be non-empty if $\pi_\a^*(L_\b)=\pi_\a^*(L_\g)$
which implies $\g=\a+\b$.

\v8\noi{\bf Note.} Since the group $\Jt$ is commutative,
we have $$|M|_\alpha \cap |M|_\beta = |M|_\alpha \cap |M|_\beta\cap
|M|_{\alpha + \beta},$$  and similarly 
in the twisted case. 
\v8
In the untwisted case, the fixed point sets have two components each, and
  we may ask about the triple intersections of the individual
  components. Our  answer was already stated in  theorem \ref{evencasei}. We will prove it in the following
  equivalent form. 

\begin{theorem}\label{6.8}  Assume that $\a$ and $\b$ are
  distinct non-zero elements of $J^{(2)}$  such
that $\lambda_2(\a,\b)=1\in\mu_2$.  Let $a,b\in J^{(4)}$ such that $2a=\a$ and
$2b=\b$. Let $\varepsilon,\mu,\nu=\pm 1$ be three signs. Then 
\begin{equation}
|M|^\varepsilon_a \cap |M|^\mu_b \cap |M|^\nu_{a + b} \neq \emptyset \
\Leftrightarrow\ \lambda_4(a,b)=\varepsilon\mu\nu\nonumber
\end{equation}\end{theorem}

\noi{\bf Note.} It follows that if $|M|^\varepsilon_a \cap |M|^\mu_b
\cap |M|^\nu_{a + b}$ is non-empty, then it is equal to all of
$|M|^\varepsilon_a \cap |M|^\mu_b$. This fact can of course be seen
directly using prop. \ref{count}. The important information in the
theorem is that it tells us when $|M|^\varepsilon_a \cap |M|^\mu_b$
intersects $|M|^+_{a + b}$ and when $|M|^-_{a + b}$.
\v8
\proof To simplify notation, we put $\g=\a+\b$ and $c=a+b$. 
Let $E$ represent a point in $|M|_\alpha \cap |M|_\beta \cap
|M|_\g$. The description of $|M|_\alpha$ in section \ref{fpv} tells us that
there exists  $\L_a\in
\pi_\a^*(L_a)\otimes \im \Phi_\a^d \subset \Pic_0(\S^\a)$ such that $E\cong\pia\L_a$;
moreover $E$ lies in $|M|_a^+$ if and only $d$ is even. (See formula
(\ref{untw2}).) Similarly we
have $E\cong\pib\L_b\cong\pig\L_c$, where $\L_b\in
\pi_\b^*(L_b)\otimes \im \Phi_\b^{d'}\subset \Pic_0(\S^\b)$ and $\L_c\in
\pi_\g^*(L_c)\otimes \im \Phi_\g^{d''}\subset \Pic_0(\S^\g)$. Thus, the theorem is
equivalent to the following lemma.

\begin{lemma}\label{6.9} One has $\lambda_4(a,b)=(-1)^{d+d'+d''}$.
\end{lemma}

The proof of this lemma will occupy the remainder of this section. 

There is a curve $\tS$ naturally double covering $\S^\a$, $\S^\b$
and $\S^{\g}$:
$$\tS = \left\{ (\xi,\eta)\in L_\a\oplus L_\b\,|\, \xi^2=1=\eta^2\right\}.$$
The projections onto the two factors induce projections
$\pi^\a \colon \tS \ra \S^\a$ and  $\pi^\b \colon \tS \ra \S^\b.$
The bilinear map $L_\a\oplus L_\b \ra L_\a\otimes L_\b$ induces the projection
$\pi^{\g} \colon \tS \ra \S^{\g}.$ 

$$\begin{array}{lll} &\tS&\\
{}^{\textstyle\pi^\a}\hskip -5pt\swarrow &\downarrow\pi^\b&\searrow^{\textstyle \pi^\g}\\
\hskip -10pt\S^\a & \S^\b & \ \ \ \ \S^\g\\
{}_{\textstyle\pi_\a}\hskip -5pt\searrow&\downarrow\pi_\b&\swarrow_{\textstyle \pi_\g}\\
&\S&
\end{array}$$

The deck-transformations of the
coverings $\pi^\a, \pi^\b, \pi^{\g}$ will be denoted respectively by $\phi^\a,
\phi^\b,\phi^\g$. Note that $\phi^\a$ (resp. $
\phi^\b$) is induced by multiplication by $-1$ in the fibers of
$L_\b$ (resp. $L_\a$), and that $\phi^\g=\phi^\a\circ
\phi^\b$. We denote the projection
$\tS\rightarrow \S$ by $\tilde \pi$, so that $$\tilde \pi=\pi_\a\circ
\pi^\a =\pi_\b\circ
\pi^\b=\pi_{\g}\circ
\pi^{\g}.$$ Notice also that the involution $\phi_\a$ of $\S^\a$ is
covered by both $\phi^\b$ and $\phi^\g$ (but of course not by
$\phi^\a$, since $\S^\a=\tS/\lag\phi^\a\rag$). Similar comments apply
to $\phi_\b$ and $\phi_\g$.

\begin{lemma}\label{6.10} One has $\pi^{\a*}(\L_a)\cong \pi^{\b*}(\L_b)\cong\pi^{\g*}(\L_c)$.
\end{lemma}
\proof Since $E\cong \pi_{\a*}(\L_a)$ lies in $|M|_\a\cap|M|_\b$, we
have $\phi_\a^*(\L_a)\cong \L_a\otimes \pi_\a^*(L_\b)$ (see the proof
of prop. \ref{trans}). Since $\pi^{\a*}\pi_\a^*(L_\b)=\tilde\pi^*(L_\b)$
is trivial, it follows that  $$\tilde
\pi^*(E)\cong\tilde
\pi^*\pia(\L_a)\cong\pi^{\a*}(\L_a\oplus \phi_\a^*(\L_a))\cong\pi^{\a*}(\L_a)\oplus\pi^{\a*}(\L_a).$$ Similarly $\tilde
\pi^*(E)\cong 
\pi^{\b*}(\L_b)\oplus\pi^{\b*}(\L_b)\cong
\pi^{\g*}(\L_c)\oplus\pi^{\g*}(\L_c)$. Since line bundles are simple,
the lemma follows.
\v8
We now turn to the computation of the Weil pairing
$\lambda_4(a,b)$. Represent $a,b\in \Jf$  by divisors
$D_a,D_b\in \Div_0(\S)$ with disjoint support, and put $D_c=D_a+D_b$. Pick $D\in Div_d(\S^\a)$ (resp. $D'\in
Div_{d'}(\S^\b)$, resp. $D''\in Div_{d''}(\S^\g)$) such that
$\pi_\a^*(D_a) +(1-\phi_\a^*)(D)$ (resp. $\pi_\b^*(D_b)
+(1-\phi_\b^*)(D')$, resp. $\pi_\g^*(D_c)
+(1-\phi_\g^*)(D'')$) represents $\L_a$ (resp. $\L_b$, resp.  $\L_c$).
 Pulling everything up to
$\tS$, we get divisors
$$F_a=\tilde\pi^{*}(D_a)+\pi^{\a*}(1-\phi_\a^*)(D)$$
$$F_b=\tilde\pi^{*}(D_b)+\pi^{\b*}(1-\phi_\b^*)(D')$$
$$F_c=\tilde\pi^{*}(D_c)+\pi^{\g*}(1-\phi_\g^*)(D'')$$
such that $F_a$ represents $\pi^{\a*}(\L_a)$, $F_b$ represents
$\pi^{\b*}(\L_b)$, and $F_c$ represents $\pi^{\g*}(\L_c)$. Since these
three bundles are isomorphic  by
lemma \ref{6.10}, there exist meromorphic functions $h_1,h_2\in
\M(\tS)$ such that $$(h_2)+F_a\,=\,F_c\,=\, (h_1)+F_b.$$ 

Let $\Nm^\a\colon \M(\tS)\ra \M(\S^\a)$ be the norm map on meromorphic
functions  associated to
the covering $\pi^\a$. The norm maps associated to the various other
coverings will similarly be denoted by
$\Nm^\b,\Nm^\g,\Nm_\a,\Nm_\b,\Nm_\g,$ and $\widetilde \Nm$.

\begin{lemma}\label{6.12} (i) Define $f,g\in \M(\S)$ by  $f=\widetilde \Nm(h_1)$,
  $g=\widetilde \Nm(h_2)$. Then $$(f)=4 \,D_a, \ \ (g)=4\, D_b.$$ 
(ii) Define $f_\a=\Nm^\a(h_1)\in \M(\S^\a)$, $f_\b=\Nm^\b(h_2)\in \M(\S^\b)$, $f_\g=\Nm^\g(h_1/h_2)\in
\M(\S^\g)$. Then $$f_\a\circ \phi_\a=-f_\a, \ f_\b\circ
\phi_\b=-f_\b, \ \ f_\g\circ \phi_\g=-f_\g.$$
\end{lemma}
\proof Using that $\phi^\a$ covers $\phi_\b$ and
$\phi_\g$, one computes  that
$$\pi^{\a*}((f_\a))=(h_1)+\phi^{\a*}(h_1)=2\,\tilde\pi^*(D_a).$$ It
follows that $(f_\a)=2\pi_\a^*(D_a)$ and hence
$(f)=(\Nm_\a(f_\a))=4D_a$, as asserted. This also shows that the
divisor $(f_\a)$ is $\phi_\a$-invariant. Therefore one has $f_\a\circ
\phi_\a=\pm f_\a$. But $f_\a$ itself cannot be $\phi_\a$-invariant,
because if it were, it would descend to a function $h\in\M(\S)$ such
that $(h)=2D_a$, which is impossible since $a$ has order
$4$. Therefore $f_\a\circ \phi_\a=-f_\a$. The other assertions  of the
lemma are proved similarly. 
\v8
By lemma \ref{6.12}(i), we can compute the Weil pairing
$\lambda_4(a,b)$ using the functions $f$ and $g$ (see the definition 
in section \ref{prel}).  Note that by lemma \ref{6.12}(ii), we have  that
$$h_1(\pi^{\a*}(1-\phi_\a^*)(D))=f_\a((1-\phi_\a^*)(D))=\frac{f_\a(D)}{(f_\a\circ
  \phi_\a )(D)}=(-1)^{\deg D}.$$

Thus  
\begin{eqnarray*}\label{calc} 
\lambda_4(a,b)&=&\frac{g(D_a)}{f(D_b)} =\frac{f(-D_b)}{g(-D_a)}=\frac{h_1(-\tilde\pi^*(D_b))}{h_2(-\tilde\pi^*(D_a))} \\
&=&\frac{h_1(-\tilde\pi^*(D_b) +(h_2))}{h_2(-\tilde\pi^*(D_a) +(h_1))}    \\
&=&\frac{h_1(\pi^{\g*}(1-\phi_\g^*)(D'') - \pi^{\a*}(1-\phi_\a^*)(D))}
       {h_2(\pi^{\g*}(1-\phi_\g^*)(D'') -
         \pi^{\b*}(1-\phi_\b^*)(D'))}\\
&=&f_\g((1-\phi_\g^*)(D'')) \frac
{f_\b((1-\phi_\b^*)(D'))}{f_\a((1-\phi_\a^*)(D))}\\
&=&(-1)^{\deg(D) +\deg (D')+\deg(D'')}=(-1)^{d+d'+d''}
\end{eqnarray*} 
where we have used Weil reciprocity in the fourth equality. This proves lemma
\ref{6.9} and hence theorem \ref{6.8}.

\section{The action of $\Jt$ on the Hecke correspondence.}\label{Heckecorr}

 We will make use of the Hecke correspondence in our analysis of the
 involutions in section \ref{Geoinvest}. This is a pair of morphisms 
\begin{center}
\begin{picture}(120,80)
\put(55,60){$\P$}
\put(0,10){$M$}
\put(100,10){$M'$}
\put(50,55){\vector(-1,-1){30}}
\put(65,55){\vector(1,-1){30}}
\put(20,45){$q$}
\put(90,45){$q'$}
\end{picture}
\end{center}
which allows one to `transfer' information from $M$ to $M'$. In this section, we describe the fixed point varieties $|\P|_\a$ of the action
of the various $\a\in\Jt$ on $\P$. 

\v8 \noi{\bf Notation.} Given a bundle $E$ over $\S$, we denote by
$E_x$ the fiber of $E$ at the point $x\in\S$. Also, for a bundle $E$
representing 
a point in the
moduli spaces $M$ or
$M'$, we use the notation
$[E]$ for that point. 
\v8 

We briefly review the construction of $\P$. (See {\em e.g.} Bertram and
Szenes \cite{BSz}.) Let
$\U$ be a Poincar\'{e} bundle over $\Sigma\times M'$. Thus, if
$[E']\in M'$, the
restriction of $\U$ to $\S\times \{[E']\}$ is isomorphic to $E'$. We can uniquely
fix $\U$ by requiring that $\det(\U|_{\{p\}\times M'})$ is an ample generator of
$\Pic(M')$.

We put $\P = \bP(\U|_{\{p\}\times M'})$ and let  $ q' : \P \ra M'$ be
the projection. Note that $q'$ is a $\bP^1$-fibration and for $[E']\in
M'$, the fiber $(q')^{-1}([E'])$ is isomorphic to the projective space
$\bP(E'_p)$. In fact, $\P$ can be viewed as the moduli space of pairs
$(E', \F)$ where $E'$ is a stable rank $2$ bundle with $\det(E')=[p]$,
  and $\F\subset E'_p$ is a one-dimensional subspace, {\em i.e.}, $\P$ is a moduli space of
semi-stable parabolic bundles. We will refer to
  $\F$ as a {\em line} in $E'_p$. Points in $\P$ will be denoted as
  $[(E',\F)]$, and we have $q'([(E',\F)])=[E']$.

The map $q$ is obtained by the operation of elementary modification at $p$. This means that we have $q([(E',\F)])=[E]$ if and
only if there is a short exact sequence (of sheaves) 
 $$0\ra E\ra E'\mapright{\lambda}\C_p\ra 0$$ such that
$\ker_p(\lambda)=\F\subset E'_p$. Here, $\C_p$ is the skyscraper sheaf
at $p$.

The group $\Jt$ acts naturally on $\P$. The action of $\a\in \Jt$ on
$\P$ sends $[(E',\F)]$ to $[(E'\otimes L_\a,\F\otimes L_\a)]$. The
morphisms $q$ and $q'$ are  $J^{(2)}$-equivariant.

Let $\a\in \Jt$ be non-zero. We now describe the fixed point variety 
$|\P|_\a$.  Recall that $\pia\colon \theta_\a^{-1}([p]) \rightarrow 
|M'|_\a$
is a double covering, where 
$\theta_\a^{-1}([p])\subset \Pic_1(\S^\a)$ consists of two `translates'
of
the Prym variety $P_\a$. 
Let $p_\alpha\in\S^\alpha$ be such that $\pi_\a(p_\alpha) =
p$. If $L\in \Pic_1(\Sigma^\alpha)$, the projection gives a canonical isomorphism
$$(L\oplus \phi_\a^*(L))_\pa \mapright \sim (\pia(L))_p.$$

\begin{proposition} \label{fixPa} We have an isomorphism $$j_{\pa}\colon
  \theta_\a^{-1}([p]) \,\mapright\sim\, |\P|_\a$$ defined 
  by $j_{\pa}(L)=[(\pia(L),(L\oplus 0)_{p_a})]$.
\end{proposition}
\proof  Put $E'=\pia(L)$. It is clear that $|\P|_\a\subset (q')^{-1}
(|M'|_\a)$. Therefore the only question is which lines in $E'_p$
correspond to fixed points of $\a$ acting on $\bP(E'_p)\cong 
(q')^{-1}([E'])$. 
Let $\psi\colon E'\mapright\sim
E'\otimes L_\a$ be the isomorphism described in remark
\ref{fundobs}. It is covered by the diagonal automorphism
$\tilde\psi=1\oplus(-1)$ of $L\oplus \phi_\a^*(L)$. Since $E'$  is stable, it is simple, 
hence any other isomorphism is a
non-zero multiple of $\psi$. Therefore $(E',\F)$ represents a point in
$|\P|_\a$ if and only if $\psi(\F)=\F\otimes L_\a$. Letting $\widetilde
\F$ denote the line in $(L\oplus \phi_\a^*(L))_\pa$ projecting down to
$\F\subset E'_p$, this condition is equivalent to $\tilde\psi(\widetilde
\F)=\widetilde
\F$. The only lines in $(L\oplus \phi_\a^*(L))_\pa$ that $\tilde\psi$ preserves are $(L\oplus 0)_\pa$
and $(0\oplus \phi_\a^*(L))_\pa$. The first line defines  the
point $j_\pa(L)$ in $\P$, and the second line defines the point
$j_\pa(\phi_\a^*(L))$. This shows that $j_\pa$ is bijective. It is clearly an
algebraic morphism,  and since its
domain is smooth, this shows $j_\pa$ is an isomorphism.

\begin{remark}\label{qeq} {\em One has $q'\circ j_\pa=\pia$. In other words,
$j_\pa$ is an isomorphism of coverings  over the identity of $|M'|_\a$. Note also that
$j_{\phi_\a(\pa)}=j_\pa \circ \phi_\a^*$.}\end{remark}
\noi{\bf Notation.} Given $a\in\Jf$ such that $2a=\a$, we denote by
$|\P|_a^+$ the component of $|\P|_\a$ containing the point
$j_\pa(\pi_\a^*(L_a)\otimes [p_\a])$. (See formula (\ref{tw}).) Note
that this point, and hence the definition of the  component
$|\P|_a^+$,
depends only on $a$, not
on the  choice of $p_\a$.
\begin{proposition}\label{q-+} One has $q(|\P|_a^+)=|M|_a^-$.
\end{proposition}
\proof Put $L=\pi_\a^*(L_a)$ and $L'=L\otimes [p_\a]$. Then
$E=\pia(L)$ represents a point in $|M|_a^+$, and $E'=\pia(L')$
represents a point in $|M'|_\a$. The short exact sequence of sheaves 
$$0\ra L\ra L'\ra \C_\pa\ra 0$$ induces the short exact sequence 
\begin{equation} 0\ra E\ra E'\,\mapright{\lambda}\,\C_p\ra 0.
\label{ses}
\end{equation}  We need to determine the
line $\F=\ker_p (\lambda)\subset E'_p$. Pulling (\ref{ses}) back to
$\S^\a$ and restricting to the fiber at $p_\a$, the map $\lambda$
becomes $$(L'\oplus \phi_\a^*(L'))_\pa\ra (\C_\pa \oplus
\phi_\a^*(\C_\pa))_\pa =\C_\pa \oplus 0.$$ This shows that $\F=\ker_p
(\lambda)$ is the projection of the line $(0\oplus
\phi_\a^*(L'))_\pa$. Hence $$[(E',
\F)]=j_\pa(\phi_\a^*(L'))=j_{\phi_\a(p_\a)}(L').$$ Since $j_\pa(L')\in
|\P|_a^+$, this shows that $[(E',
\F)]$ lies
in $|\P|_a^-$ (see  remark \ref{qeq}). Recalling that $q([(E', \F)]=[E]$, it follows that 
$q(|\P|_a^-)=|M|_a^+ $,   and also that  $q(|\P|_a^+)=|M|_a^-$.
\v8

The following observation will be used in section \ref{tracecomp}.
\begin{proposition}\label{qtriv} Let $\nu$ be the relative cotangent
  sheaf of $q'\colon \P\ra M'$. Then the restriction of $\nu$ to
  $|\P|_\a$ is numerically trivial.  
\end{proposition}

\proof Note that $\nu|_{|\P|_\a}$ is the dual of the normal bundle,
$N$, say, of
the inclusion $|\P|_\a\subset (q')^{-1}(|M'|_\a)$. By
prop. \ref{fixPa}, it suffices to show
that $j_\pa^*(N) $ is numerically trivial. Let $\Lambda$ be a
Poincar\'e  bundle over 
$\Pic_1(\S^\a)\times \S^\a$. For $x\in \S^\a$, let $\Lambda_{x}$ denote its
restriction to $\theta_\a^{-1}([p])\times \{x\}$. We
have a commutative diagram
$$\begin{array}{cccc}
\bP( \Lambda_{p_\alpha} \oplus \Lambda_{\phi_\a(\pa)} )
  &\mapright{\Pi_{\a*}} & (q')^{-1}(|M'|_\a)&\ \subset \ \P \\ 
\downarrow & & \downarrow \, q' &\\
\theta_\a^{-1}([p]) & \mapright{\pia} & |M'|_\a&
\end{array}$$ Here, $\Pi_{\a*}$ is the obvious map covering $\pia$. (A
point in $\bP( \Lambda_{p_\alpha} \oplus \Lambda_{\phi_\a(\pa)} )$ is
a point $[L]\in \theta_\a^{-1}([p])$ together with a line in $
L_{\pa}\oplus L_{\phi_\a(\pa)}=(L\oplus \phi_\a^*(L))_\pa$. This is
sent by 
$\Pi_{\a*}$ to the point represented by $\pia(L)$ and the induced line
in $(\pia(L))_p$.)  

Let $s_{\pa}$ be the section of the fibration on
the left defined by $s_\pa(L)=L_{\pa}\oplus 0$, for $[L]\in
\theta_\a^{-1}([p])$. Then $j_\pa=\Pi_{\a*}\circ s_{\pa}$, and hence
the inclusion $|\P|_\a\subset (q')^{-1}(|M'|_\a)$ corresponds to the
inclusion of the image of $s_\pa$ in $\bP( \Lambda_{p_\alpha} \oplus
\Lambda_{\phi_\a(\pa)} )$. This shows 
$$j_\pa^*(N) \cong \Lambda_\pa^*\otimes \Lambda_{\phi_\a(\pa)}.$$
Tensoring $\Lambda$ by the pull-back of a bundle over $\Pic_1(\S^\a)$
if necessary, we may assume $\Lambda_\pa$ is trivial. Hence
$j_\pa^*(N)$ is numerically 
trivial, proving the proposition.

\section{Investigation of the involutions $\rho_a$ and $\rho'_\a$.}
\label{Geoinvest}

Let $a\in\Jf$ such that  $2a=\a\neq0$. Recall that the involution
$\rho_a$ is the lift of $\a$ to $\L$ which acts as the identity  over
the fixed point component $|M|_a^+$.

\begin{proposition}\label{rhominus} The involution $\rho_a$  acts as minus the identity  over
the fixed point component $|M|_a^-$.
\end{proposition} 
\proof Let
$E'$ represent a point $[E']$ in $|M'|_\alpha$. Now
$\alpha$ acts on the fiber of $q'$ over $[E']$ and from our description of
$|\P|_\alpha$, we have that
$\alpha$ has exactly two fixed points on $(q')^{-1}([E'])$. Consider now
$q^*\L|_{(q')^{-1}([E'])}$ with its lift of $\alpha$ induced by
$\rho_a$. Let $s_1,s_2$ be the signs by which $\rho_a$ acts  over the
two fixed points. By lemma 2.1 in
\cite{BSz} we have that
$q^*\L|_{(q')^{-1}([E'])}$ is isomorphic to $\cO(1)$ over
$(q')^{-1}([E'])\cong \bP(E'_p)$. From this we conclude that $s_1s_2 = -1$,
proving the proposition.

\begin{remark}\label{coc} {\em  At this point, it follows easily  that the
    alternating form  of the extension $\E$ generated by the involutions $\rho_a$
    is equal to the order $2$ Weil pairing $\lambda_2$. Indeed, recall
    that the alternating form  is defined by the commutator pairing
$c(\a,\b)=\rho_a\rho_b\rho_a^{-1}\rho_b^{-1}$ where $a$ is  a lift of
$\alpha$ and $b$ is  a
lift of $\beta$. Assume first that $\lambda_2(\alpha, \beta)=1$.  
Let us then evaluate $\rho_a\rho_b\rho_a^{-1}\rho_b^{-1}$ in a point in
$|M|^+_a$. Recalling from \ref{4.6}(ii) that $\beta$ preserves this component of $|M|_\alpha$,
 we get that
$$c(\alpha,\beta) = 1 \,\rho_b \,(1)^{-1}\, \rho_b^{-1} = 1.$$ 
If however $\lambda_2(\alpha, \beta)=-1$, then $\beta$ exchanges the two
components, and  $$c(\alpha,\beta) = 1 \,\rho_b \,(-1)^{-1}\, \rho_b^{-1} =
-1.$$
Thus $c=\lambda_2$, as asserted.}
\end{remark}

\begin{theorem}\label{8.4} We have
  $\rho'_\alpha
\rho'_\beta=\lambda_2(\a,\b)\rho'_{\alpha+\beta}$. 
 \end{theorem} 
\proof We may assume none of the classes $\a$, $\b$, $\a+\b$,  is zero, the
result being obvious otherwise. Consider first the case
$\lambda_2(\a,\b)=-1$. Then by prop. \ref{6.5} we have that the triple
intersection $
|M'|_\alpha \cap |M'|_\beta \cap |M'|_{\alpha + \beta}$ is non-empty.
Since by definition $\rho'_\alpha$ acts as minus the identity on the
fiber over $|M'|_\alpha$,  it
follows that  $$\rho'_\alpha \rho'_\beta=-\rho'_{\alpha +
  \beta}, $$ proving the result in this case.  

Now consider the case
$\lambda_2(\a,\b)=1$. Pick $a,b\in \Jf$ such that $2a=\a$ and
$2b=\b$, and consider the involutions $\rho_a^{\otimes 2}$, 
$\rho_b^{\otimes 2}$, and $\rho_{a+b}^{\otimes 2}$,  acting on
$\L^2$. In fact, those involutions depend only on $\a$ and $\b$, and
not on the choice of $a$ and $b$. By prop. \ref{6.5} we have that the triple
intersection $
|M|_\alpha \cap |M|_\beta \cap |M|_{\alpha + \beta}$ is
non-empty. Note that $\rho_a^{\otimes 2}$ acts as the identity over
both components of $|M|_\a$.  Hence 
\begin{equation}\label{rhoH} \rho_a^{\otimes 2} \rho_b^{\otimes 2}=\rho_{a+b}^{\otimes 2}.
\end{equation} Now consider the Hecke correspondence.  From Corollary
2.2 in \cite{BSz} (see also Lemma 10.3 in \cite{BLS}) we have that the
canonical bundle $K_{\P}$ of $\P$ satisfies 
\begin{equation}\label{Hecke}
K_{\P} \cong (q')^*({\L'}^{-1})\otimes q^*(\L^{-2}).
\end{equation}
  From Proposition
\ref{fixPa} we see that $|\P|_\alpha$ has odd codimension, hence $\alpha$ acts by
$-1$ on the restriction of
$K_\P$ to
$|\P|_\alpha$. Our lifts  $\rho_a^{\otimes 2}$ and $\rho'_\alpha$ thus
make the isomorphism (\ref{Hecke}) a $J^{(2)}$-equivariant
isomorphism. The action of $\Jt$ on 
$K_{\P}$ is obviously a group action. This enables us to compare the lift $\rho'_\alpha$ acting on ${\L'}$
and the lift $\rho_a^{\otimes 2}$ acting on $\L^{2}$. Thus
(\ref{rhoH}) 
implies $$\rho'_\alpha \rho'_\beta=\rho'_{\alpha +
  \beta},$$ proving the result in the case $\lambda_2(\a,\b)=1$. This
 completes the proof.
\v8\noi{\bf Note.} The equivariance of the isomorphism (\ref{Hecke})
is the reason why we defined $\rho'_\alpha$ to be the lift which acts
as {\em minus} the identity over the fixed point set.

\begin{theorem}\label{8.5} We have $ \rho_a \rho_b=\lambda_4(a,b)\rho_{a+b}$.
\end{theorem}

\proof  We first deal with the case where $\lambda_4(a,b)=\pm 1$, or,
equivalently, 
$\lambda_2(\a,\b)=1$. If $\a=\b=0$, there is nothing to show. If
$\a=\b\neq 0$, then $\rho_a$ and $\rho_b$ are lifts of the same class,
hence  $\rho_a=\pm \rho_b$. By  prop. \ref{rhominus}, we have $\rho_a=\rho_b$ if and only if $a$
and $b$ define the same component of $|M|_\a$, which in turn is
equivalent, by prop. \ref{4.6}(i), to $\lambda_2(b-a,\a)=1$.  But if $\a=\b$ then
$\lambda_2(b-a,\a)=\lambda_4(a,b)$ and  $\rho_{a+b}$ is the
  identity (by definition). This proves the result in the case $\a=\b$.   Finally, if
  $\a,\b,$ and $\a+\b$ are all three non-zero, the triple intersection
   $ |M|_\alpha \cap
|M|_\beta \cap |M|_{\alpha + \beta}$ is non-empty, and we can compute $\rho_a \rho_b 
\rho_{a+b}$ in the fiber over an intersection point. By theorem
\ref{evencasei} and prop. \ref{rhominus}, it follows that $\rho_a \rho_b 
\rho_{a+b}=\lambda_4(a,b)$, completing the proof in the case
$\lambda_4(a,b)=\pm 1$.

The remainder of this section is devoted to the proof in the case  where  $\lambda_4(a,b)=\pm
i.$ We
will again use the notations $\g=\a+\b$ and $c=a+b$.

Let the bundle 
$E'$ represent a point $[E']$ in $|M'|_\alpha\cap |M'|_\b\cap
|M'|_\g$. The three involutions $\a,\b,\g$ induce involutions on
$(q')^{-1}([E'])\subset \P$; recall that $(q')^{-1}([E'])$  is
identified with  the
projective space $\bP(E'_p)$. Each of these involutions has two fixed
points on $\bP(E'_p)$; these are precisely the intersection points of
$\bP(E'_p)$ with
the fixed point varieties $|\P|_\a,|\P|_\b$, and $|\P|_\g$.  Note that 
$\rho_a$ acts as $\mp 1$ on the fiber of $q^*\L$ at the intersection point
of $\bP(E'_p)$ with 
the component $|\P|_a^{\pm}$, since
$q(|\P|_a^+)=|M|_a^-$ by prop. \ref{q-+}. Let
$\F_a^{\pm}$,  $\F_b^{\pm}$, $\F_c^{\pm}$ be the lines in $E'_p$
corresponding to the intersection points of $\bP(E'_p)$ with  the components $|\P|_a^{\pm},|\P|_b^{\pm}$,
and $|\P|_c^{\pm}$. As already used in the proof
of prop. \ref{rhominus}, the restriction of $q^*\L$ to $\bP(E'_p)$ is
the bundle $\cO(1)$.  It will be convenient to transfer the calculation
to the tautological bundle $\cO(-1)$, whose fiber over a point
represented by a line $\F$ is that line. Of course, $\cO(-1)$ is the restriction of
$q^*\L^{-1}$ to $\bP(E'_p)$. For $a\in\Jf$, let us denote
by 
$\hat\rho_a$ the involution  
$\rho_a^{\otimes(-1)}$ acting on $\L^{-1}$. Then  $\hat\rho_a$ acts as $\mp 1$ on the line
$\F_a^{\pm}$, and similarly for $\hat\rho_b$ and $\hat\rho_c$ on the lines
$\F_b^{\pm}$ and $\F_c^{\pm}$.

 It follows easily from this description (or from
the computation of the alternating form in remark \ref{coc}) that one has
$\hat\rho_a\hat\rho_b=\varepsilon \hat\rho_c$ where $\varepsilon\in\{\pm i\}$. (In
fact, the involutions $\hat\rho_a,
\hat\rho_b$ generate a quaternion subgroup $Q_8\subset Sl_2(\C)$,
covering the commutative subgroup $\Z/2\times \Z/2 \subset PSl_2(\C)$
generated by $\a,\b$.)  Of
course, the sign of $\varepsilon$ is determined by the relative position of the
six lines. The following lemma computes this relative position in terms of the Weil pairing
$\lambda_4(a,b)$.

\begin{lemma}\label{ml} Let $\lambda=\lambda_4(b,a)\in\{\pm i\}$. There is an isomorphism of
  $E'_p$ with $\C^2$ sending the six lines $\F_a^{+}$, $\F_a^{-}$,
  $\F_b^{+}$, $\F_b^{-}$, $\F_c^{+}$, $\F_c^{-}$, to the lines
  generated by the vectors 
$$\left(\begin{array}{c} 1\\0 \end{array}\right), \ 
\left(\begin{array}{c} 0\\1 \end{array}\right), \ 
\left(\begin{array}{c} 1\\1 \end{array}\right), \ 
\left(\begin{array}{c} 1\\-1\end{array}\right), \ 
\left(\begin{array}{c} 1\\-\lambda \end{array}\right), \ 
\left(\begin{array}{c} 1\\ \lambda \end{array}\right) \ . 
$$
\end{lemma}

The proof will be given later. Assuming lemma \ref{ml} for the moment,
we see that $\hat\rho_a,
\hat\rho_b$ and $\hat\rho_c$ correspond to the matrices 
$$T_a=\left( \begin{array}{cc} -1& 0 \\ 0 & 1 \end{array}\right), \ \ 
T_b=\left( \begin{array}{cc} 0& -1 \\ -1 & 0 \end{array}\right), \ \
T_c=\left( \begin{array}{cc} 0& -\lambda \\ \lambda & 0
  \end{array}\right).$$

Note that $T_aT_b=\lambda T_c$, and hence
$\hat\rho_a\hat\rho_b=\lambda \hat
\rho_c$ and $\rho_a \rho_b=\lambda^{-1}\rho_c$.  Thus the remaining case of theorem \ref{8.5} follows directly
from lemma \ref{ml}.
\v8
Now let us prove  lemma \ref{ml}. The proof uses again  the coverings $\S^\a$,
$\S^\b$, $\S^\g$, and their common covering $\tS$ (see section \ref{triple}). 

Choose a point $\tilde p\in\tilde\pi^{-1}(p)\subset \tS$ and put
\begin{equation}
p_\a=\pi^\a(\tilde p), \ p_\b=\pi^\b(\tilde p), \ p_\g=\pi^\g(\tilde
p).
\label{palpha}
\end{equation}  
Since  $[E']\in |M'|_\alpha\cap |M'|_\b\cap
|M'|_\g$, there exist line bundles $\L_a$ over $\S^\a$, $\L_b$ over
$\S^\b$,   and $\L_c$ over
 $\S^\g$, such that $
E'\cong
\pia(\L_a)\cong
\pib(\L_b)\cong
\pig(\L_c)$. We can fix $\L_a$ (resp. $\L_b$, resp. $\L_c$) uniquely up to isomorphism by requiring
that $\L_a\in 
\pi_\a^*(L_a)\otimes[p_\a]\otimes \im\Phi_\a^0$ (resp. $\L_b\in 
\pi_\b^*(L_b)\otimes[p_\b]\otimes \im\Phi_\b^0$, resp.   $\L_c\in 
\pi_\g^*(L_c)\otimes[p_\g]\otimes \im\Phi_\g^0$) (see formula (\ref{tw})).

Let us denote the bundle $\pi^{\g*}(\L_c)$ on $\tS$ by $L$. Being a
pull-back bundle, $L$ has a canonical involution, $C$, say, covering the involution $\phi^\g$ on
$\tS$, and such that $L/\lag C\rag$  is the bundle $\L_c$ on
$\S^\g=\tS/\lag\phi^\g\rag$.  Proceeding as in lemma \ref{6.10}, it is easy to check that the 
bundles $\pi^{\a*}(\L_a)$ and $\pi^{\b*}(\L_b)$ are isomorphic to $L=\pi^{\g*}(\L_c)$. Therefore $L$ also has canonical involutions $A$
and $B$, covering $\phi^\a$ and $\phi^\b$, respectively, such that
$L/\lag A\rag\cong \L_a$ and $L/\lag B\rag\cong \L_b$. 

\begin{lemma}\label{ABC}  One has $AB=\lambda_4(b,a)C$.
\end{lemma} 
\proof As in the proof of lemma \ref{6.9}, represent
$\pi^{\a*}(\L_a)$, $\pi^{\b*}(\L_b)$, and $L=\pi^{\g*}(\L_c)$, by
divisors $F_a$, $F_b$, $F_c$, respectively, such that 
$$F_a=\tilde\pi^{*}(D_a)+\pi^{\a*}(p_\a+(1-\phi_\a^*)(D))$$
$$F_b=\tilde\pi^{*}(D_b)+\pi^{\b*}(p_\b+(1-\phi_\b^*)(D'))$$
$$F_c=\tilde\pi^{*}(D_c)+\pi^{\g*}(p_\g+(1-\phi_\g^*)(D''))$$
where $D_a,D_b\in \Div_0(\S)$ represent $a,b\in \Jf$, $D_c=D_a+D_b$, $D\in Div_0(\S^\a)$, $D'\in 
Div_{0}(\S^\b)$, and $D''\in Div_{0}(\S^\g)$. As before, since the three
bundles are isomorphic, there exist meromorphic functions $h_1,h_2\in
\M(\tS)$ such that $$(h_2)+F_a\,=\,F_c\,=\, (h_1)+F_b.$$ 

The action of our involutions $A,B$ and $C$ on $L$ can be
described on local sections as
follows. Since $L=\cO(F_c)$, a local section over some open set $U\subset
\tS$ is just a meromorphic function $s$ on $U$ such that $(s)+F_c|_U\geq 0$. The action of $C$ is
simply given by $$ C\ : s \mapsto s\circ \phi^\g,$$ since the divisor
$F_c$ was pulled back from $\S^\g$, and $C$ is the canonical
involution of the pull-back bundle. 

The involution $A$ is nothing but the
canonical involution of the pull-back bundle $\pi^{\a*}(\L_a)=\cO(F_a)$,
conjugated by an 
isomorphism with $L=\cO(F_c)$.  Since $(h_2)=F_c-F_a$, multiplication
by $h_2$ gives such  
an isomorphism $\cO(F_c)\,\mapright\sim \,\cO(F_a)$. Therefore the action
of $A$ on local sections of $L=\cO(F_c)$ is
$$ A\ : s \mapsto ((s h_2)\circ \phi^\a)h_2^{-1}=(s\circ
\phi^\a)\,k_A= (s\
k_A^{-1})\circ \phi^\a$$
where we have put $k_A=(h_2\circ \phi^\a)/h_2$. (N.b., one may think about this as
follows: $s\circ
\phi^\a$ is a local section of $\phi^{\a*}L=\cO(\phi^{\a*}F_c)$, and
multiplication by $k_A$ describes an isomorphism
$\cO(\phi^{\a*}F_c)\,\mapright\sim \,\cO(F_c)$.)
 
Similarly, $B$ acts on local sections of
$\cO(F_c)$ as $$ B\ : s \mapsto ((s h_1)\circ \phi^\b)h_1^{-1}=(s\circ
\phi^\b)\,k_B= (s\
k_B^{-1})\circ \phi^\b$$
where $k_B=(h_1\circ \phi^\b)/h_1$.

Put \begin{equation} \lambda= \frac{k_B}{k_A} =
\frac{h_2}{h_1} \,\frac {h_1\circ \phi^\b}{h_2\circ
  \phi^\a}\label{la1}
\end{equation} Note that $\lambda$ is a constant, since
$(k_B)=(k_A)$. Since  $AB$ acts on
local sections by 
$$AB\ : s\mapsto (((s\circ
\phi^\b)\,k_B) \circ \phi^\a)\, k_A=((s\circ
\phi^\b)\,k_B\,k_A^{-1})\circ \phi^\a =\lambda \, s\circ \phi^\g,$$
we have  $AB=\lambda C$.

Now let us show that $\lambda=\lambda_4(b,a)$. As in section
\ref{triple}, 
we use the functions $f,g\in \M(\S)$ defined by  $f=\widetilde
  \Nm(h_1)$ and $g=\widetilde \Nm(h_2)$. A computation shows that the
  statements of lemma \ref{6.12} hold word for word. We can therefore
  compute $\lambda_4(a,b)=g(D_a)/f(D_b)$ exactly as before. Note that
  the divisors $D,D',D''$ have degree zero, so that  the terms 
  involving 
the functions $f_\a,f_\b, f_\g$  are now equal to $1$. We thus obtain 

\begin{eqnarray*}
\lambda_4(a,b)&=&\frac{h_1(\pi^{\g*}(p_\g)-\pi^{\a*}(p_\a))}
{h_2(\pi^{\g*}(p_\g)-\pi^{\b*}(p_\b))} 
=\frac{h_1(\phi^\g(\tilde p)-\phi^\a(\tilde p))}{h_2(\phi^\g(\tilde p)-\phi^\b(\tilde
  p))} \\
&=& 
\frac{h_1}{h_2} \,\frac {h_2\circ \phi^\a}{h_1\circ
  \phi^\b} \,(\phi^\g(\tilde p)) 
\end{eqnarray*} 
where we have used (\ref{palpha}) in the last but one step.
Comparing this with formula (\ref{la1}), we have
$\lambda=\lambda_4(a,b)^{-1}=\lambda_4(b,a)$, proving lemma \ref{ABC}.

We return to the proof of lemma \ref{ml}. Recall that 
\begin{equation}
E'\cong
\pia(\L_a)=(\L_a\oplus\phi_\a^*\L_a)/\lag\phi_\a\rag. \label{uniqueflag}
 \end{equation} 
Since $\L_a\in 
\pi_\a^*(L_a)\otimes[p_\a]\otimes \im\Phi_\a^0$, we have
$$j_\pa(\L_a)\in |\P|_a^+,$$ where $j_\pa$ is the isomorphism of 
prop. \ref{fixPa}. Thus,  the line $\F_a^+\subset
E'_p$ (representing the unique intersection point in $\bP(E'_p)\cap
|\P|_a^+$)  is the image of the line  $(\L_a\oplus 0)_{\pa}\subset 
(\L_a\oplus\phi_\a^*\L_a)_{\pa}$ under the
natural projection. Also, the line $\F_a^-$ is the image of $(0\oplus
\phi_\a^*(\L_a))_{\pa}$. Note that since $E'$ is stable, the isomorphism
(\ref{uniqueflag}) is unique up to scalar multiples, hence the lines
$\F_a^{\pm}$ are  
well-determined by this description. 

Similarly, using the isomorphism of $E'$ with $\pib(\L_b)$ and with
$\pig(\L_c)$, the lines  $\F_b^+$ and $\F_b^-$ correspond to the projections of
the lines $(\L_b\oplus 0)_{\pb}$ and $(0\oplus
\phi_\b^*(\L_b))_{\pb}$, and the lines  $\F_c^+$ and $\F_c^-$ correspond to the projections of
the lines $(\L_c\oplus 0)_{\pg}$ and $(0\oplus
\phi_\g^*(\L_c))_{\pg}$.

Let us now understand the relative position of the six lines in
$E'_p$. 

Note that since $AB=-BA$, the involution  $B$ of $L$ covering the involution
$\phi^\b$ on $\tS$ induces a map
$$L/\lag -A\rag\,\mapright{}\,L/\lag A\rag$$ covering the involution $\phi_\a$ on
$\S^\a$. We may choose isomorphisms $\L_a\cong L/\lag A\rag$ and
$\phi_\a^*\L_a\cong L/\lag -A\rag$ such that this map becomes  the canonical
map $\phi_\a^*\L_a \ra \L_a$ covering $\phi_\a$.  Therefore
$\pia(\L_a)$ is isomorphic to the bundle $E_a$ defined by 

$$E_a\ = \ L\oplus L {\Big/} \bigl< 
\left( \begin{array}{cc} A& 0 \\ 0 & -A \end{array}\right), 
\left( \begin{array}{cc} 0& B \\ B & 0 \end{array}\right)\bigr>.$$

(Notice that the two matrices commute, so that $E_a$ is indeed a
well-defined bundle on $\S=\tS/\langle\phi^\a,\phi^\b\rag$.) Moreover, the lines
$\F_a^+$ and $\F_a^-$ in $E'_p$ correspond, {\em via} an isomorphism
$E'\cong E_a$ (which is unique up to scalar multiples), to the
images of the lines $(L\oplus 0)_{\tilde p}$
and $(0\oplus L)_{\tilde p}$ under the projection from
$L\oplus L$ to $E_a$. Here, we have used that $\pi^\a(\tilde p)=\pa$
by our choice of $p_\a$ in (\ref{palpha}).  

Similarly, $\pib(\L_b)$ is isomorphic to the bundle $E_b$ defined by 

$$E_b\ = \ L\oplus L {\Big/} \bigl< 
\left( \begin{array}{cc} B& 0 \\ 0 & -B \end{array}\right), 
\left( \begin{array}{cc} 0& A \\ A & 0 \end{array}\right)\bigr>,$$ and
since $\pi^\b(\tilde p)=\pb$,
the lines $\F_b^+$ and $\F_b^-$ in $E'_p$ correspond to the
images of the same lines $(L\oplus 0)_{\tilde p}$
and $(0\oplus L)_{\tilde p}$, but  now projected from
$L\oplus L$ to $E_b$.

Lastly, $\pig(\L_c)$ is isomorphic to the bundle $E_c$ defined by 
$$E_c\ = \ L\oplus L {\Big/} \bigl< 
\left( \begin{array}{cc} C& 0 \\ 0 & -C \end{array}\right), 
\left( \begin{array}{cc} 0& A \\ A & 0 \end{array}\right)\bigr>, $$
and as before, since $\pi^\g(\tilde p)=\pg$, the lines $\F_c^+$ and $\F_c^-$ in $E'_p$ correspond to the
images of the lines $(L\oplus 0)_{\tilde p}$
and $(0\oplus L)_{\tilde p}$ under the projection from
$L\oplus L$ to $E_c$.

The three bundles $E_a$,$E_b$, $E_c$ are all isomorphic to
$E'$.  In fact, we have isomorphisms
$\psi_{X}\colon E_b\mapright\sim  E_a$ and  $\psi_{Y}\colon
E_c\mapright\sim 
E_a$ induced by the endomorphisms of $L\oplus L$ defined   by the
matrices $$X=\left( \begin{array}{cc} 1& 1 \\ 1 & -1
  \end{array}\right), \ \ \ Y=\left( \begin{array}{cc} 1& 1 \\
    -\lambda & \lambda \end{array}\right).$$ 
The verification, which uses that $AB=\lambda C$ by lemma \ref{ABC},
is left to the reader.

Let us identify the fiber $(L\oplus
L)_{\tilde p}$ in the obvious way with $\C\oplus \C$ and consider the
isomorphism $$E'_p\,\mapright\sim\, (E_a)_p
\stackrel{\sim}\longleftarrow (L\oplus
L)_{\tilde p} = \C\oplus \C,$$ where the first map is induced by an
isomorphism $E'\cong E_a$ and the second map is the projection. 
Then the six lines $\F_a^{+}$, $\F_a^{-}$,
  $\F_b^{+}$, $\F_b^{-}$, $\F_c^{+}$, $\F_c^{-}$, correspond to the lines
in $\C\oplus \C$  generated by the vectors 
$$\left(\begin{array}{c} 1\\0 \end{array}\right), \ 
\left(\begin{array}{c} 0\\1 \end{array}\right), \ 
X\left(\begin{array}{c} 1\\0 \end{array}\right) , \ 
X\left(\begin{array}{c} 0\\1 \end{array}\right), \ 
Y\left(\begin{array}{c} 1\\0 \end{array}\right) , \ 
Y\left(\begin{array}{c} 0\\1 \end{array}\right) \ . 
$$
These vectors are precisely the ones in the statement of lemma
\ref{ml}. 
This proves lemma \ref{ml}, and completes the proof of theorem
\ref{8.5}.

\section{The trace computation.}\label{tracecomp}

In this section, we prove theorem \ref{1.2} by computing the trace of ${\rho'_\alpha}^{\otimes k/2}$
and $\rho_a^{\otimes k}$ using the Lefschetz-Riemann-Roch fixed point
formula. In the twisted case, the computation is rather straightforward,
since $M'$ is smooth, and the relevant cohomology classes 
on the fixed
point set are given in \cite{NR}.  
This computation has been done in a different context by Pantev \cite{Pa}. We
repeat the calculation here and in the process we correct a misprint
in his formula. In the untwisted case, the moduli space and the fixed point
sets are not smooth. We circumvent this problem by transferring the
computation to $\P$, using some results of \cite{BSz}.

\v8\noi{\bf Note.} Beauville \cite{Be3} has recently 
computed the traces in the untwisted case in a different way by transferring the calculation to
$M'$. \footnote{Our computation was  
done independently of his.} Beauville 
considers, more generally,  rank $r$ bundles, and his formula agrees with ours in the
case $r=2$. He does not, however,  choose lifts to the line bundle
$\L$, and his result (for $r=2$)   concerns only the case
$k\equiv 0$ mod $4$, where one has a group action of $\Jt$. 

\begin{proposition}
The trace of the involution ${\rho'_\alpha}^{\otimes k/2}$ is given by
$$\tr( {\rho'_\alpha}^{\otimes k/2})= (-1)^{k/2}\left( \frac{k+2}{2}\right)^{g-1}.$$
\end{proposition}

\proof The Lefschetz-Riemann-Roch fixed point formula \cite{AS} states that
\[ \tr({\rho'_\alpha}^{\otimes k/2}) = \ch({\L'}|_{|M'|_\alpha})^{k/2}\,
      \ch(\lambda_{-1}N_\alpha)^{-1}\,\td(|M'|_{\alpha})\cap [|M'|_\alpha].\]
Here $\td(|M'|_{\alpha})$ is the Todd class, 
$N_\alpha$ is the conormal bundle of $|M'|_\alpha$, $\lambda_{t}$ is
the operation defined by $\lambda_t E= \sum
t^i \Lambda^i E$, and $$\ch(E)=Ch(E_+ -E_-),$$ where $Ch$ is the Chern
character, and,  for any
$\Z/2$-equivariant bundle $E$ over the fixed point set $|M'|_\alpha$,
$E_+$ and $E_-$ are the  $\pm 1$-eigenbundles.

Since the fixed point set $|M'|_\alpha$ is isomorphic to the Prym
variety $P_\alpha$, its Todd class is $1$.
The cohomology class $\ch(\lambda_{-1}N_\alpha)$ was computed in Proposition
4.2 in \cite{NR}. Let $\Theta$ be the restriction to $P_\alpha$ of the
principal polarization of $J_0(\S^\alpha)$. The result of Narasimhan and Ramanan
then states that
\begin{equation}
\ch(\lambda_{-1}N_\alpha) = 2^{2(g-1)}e^{-2\Theta}. \label{Thet2}
\end{equation}
By the construction of ${\L'}$ in \cite{DN}, we have the following relation
between
$\Theta$ and
$Ch({\L'}|_{|M'|_\alpha})$. 
\begin{equation}
Ch({\L'}|_{|M'|_\alpha}) = e^{2\Theta}. \label{Thet1}
\end{equation} Note that $\ch({\L'}|_{|M'|_\alpha})=-
Ch({\L'}|_{|M'|_\alpha})$,  by our definition of $\rho'_\a$. 
Hence we get that
\[ \tr({\rho'_\alpha}^{\otimes k/2})= (-1)^{k/2}2^{-2(g-1)} e^{(k+2)\Theta}\cap
[P_\alpha].\] Using Corollary 4.16 in \cite{NR}, which states that
\begin{equation}
\Theta^{g-1}\cap [P_\alpha] = (g-1)!2^{g-1}, \label{Thet3}
\end{equation}
the result follows.

\begin{proposition}
The trace of the involution $\rho_a^{\otimes k}$ is given by
$$\tr(\rho_a^{\otimes k}) = \frac{1+(-1)^k}{2}\left(
\frac{k+2}{2}\right)^{g-1}.$$
\end{proposition}

\proof
The morphism $q$ induces a $\Jt$-equivariant morphism
$$q^* : H^0(M,\L^k) \ra H^0(\P,q^*\L^k).$$
According to \cite{BSz} we have that this morphism is an isomorphism and that
$$H^i(\P,q^*\L^k) = 0,$$
for $i>0$. Hence we just need to apply the  Lefschetz-Riemann-Roch
fixed point theorem to $(\P,q^*\L^k)$. The fixed point set $|\P|_\a$
has two components, $|\P|_a^+$ and $|\P|_a^-$, each of which is
isomorphic to the Prym variety $P_\a$, and hence has trivial Todd
class. In order to understand $q^*\L|_{|\P|_\a}$ and the conormal
bundle $N(|\P|_\alpha)$, consider the following exact sequence
$$ 0 \ra (q')^*T^*_{M'} \ra T^*_\P \ra \nu \ra 0,$$
where $\nu$ is the relative cotangent sheaf of $q'$. We 
conclude that
$$K_\P \cong \nu \otimes (q')^*K_{M'}\cong \nu \otimes (q')^*
(\L')^{-2},$$ since $K_{M'} \cong (\L')^{-2}$ \cite{DN}. 
By cor. 2.2 in \cite{BSz} (see equation (\ref{Hecke}) in the proof of
theorem \ref{8.4}), we get that
$$q^*\L^2 \cong (q')^*{\L'} \otimes \nu^{-1}.$$
Now recall from  prop. \ref{qtriv} that the line bundle $\nu|_{|\P|_\a}$ is 
numerically trivial. It follows that  
$$Ch(q^*\L|_{|\P|_a^{\pm}})=\bigl(Ch((q')^*{\L'}|_{|\P|_a^{\pm}})\bigr)^{1/2}=e^\Theta$$
  where in the last equality we have used formula (\ref{Thet1}), after
  identifying both $|\P|_a^+$ and
  $|\P|_a^-$ with the Prym variety $P_\a$. Next, observe that $\alpha$
  acts as $-1$ on 
  $\nu|_{|\P|_\a}$, since one has an exact
  sequence (of bundles over $|\P|_\a$)
\begin{equation}
0\ra (q')^*N_\a \ra N(|\P|_\alpha)\ra \nu|_{|\P|_\alpha}\ra
0.\label{exse}
\end{equation} where $N_\a=N(|M'|_\alpha)$ as before. Therefore 
$$\ch(\lambda_{-1}\nu|_{|\P|_a^{\pm}})=\ch(1-\nu|_{|\P|_a^{\pm}})=Ch(1+1)=2,$$
    and 
    the exact sequence (\ref{exse}) 
 gives us that  $$\ch(\lambda_{-1}(N(|\P|_a^{\pm})) )=
2\, \ch(\lambda_{-1}(N_\a))=2\cdot
2^{2(g-1)}e^{-2\Theta}, $$ where we have used formula (\ref{Thet2}) in
the last equality, after again identifying $|\P|_a^+$ and
  $|\P|_a^-$ with $P_\a$.

Now recall that $\rho_a$ acts with opposite signs on the restriction
of $q^*\L$ to the two components $|\P|_a^+$ and
  $|\P|_a^-$. In fact, it acts as $\mp 1$ over $|\P|_a^{\pm
    }$, by prop. \ref{q-+}. Therefore $$\ch(q^*\L|_{|\P|_a^{\pm }})=\mp
  Ch(q^*\L|_{|\P|_a^{\pm }})=\mp e^\Theta.$$ Putting everything together,  the fixed point formula gives  
$$\tr(\rho_a^{\otimes k}) = (1 + (-1)^k) e^{k\Theta} \, 
2^{-1} 2^{-2(g-1)} e^{2\Theta} \cap
[P_\alpha].$$
The proposition now follows as in the twisted case from formula
(\ref{Thet3}).

\end{document}